\documentclass[11pt]{article}
\usepackage{amsmath, amsfonts, amssymb}
\usepackage{graphicx, amsthm,color}

\newtheorem{theorem}{Theorem}[section]
\newtheorem{lemma}[theorem]{Lemma}

\newtheorem{proposition}[theorem]{Proposition}

\theoremstyle{definition}

\theoremstyle{remark}
\newtheorem{remark}[theorem]{Remark}
\numberwithin{equation}{section}

\newcommand{\R}{\mathbb{R}}

\providecommand{\ri}{\mathop{\rm ri}\nolimits}
\providecommand{\into}{\mathop{\rm int}\nolimits}
\providecommand{\dom}{\mathop{\rm dom}\nolimits}

\providecommand{\epi}{\mathop{\rm epi}\nolimits}
\providecommand{\argmin}{\mathop{\rm argmin}}

\providecommand{\Liminf}{\mathop{\rm Liminf}}
\providecommand{\Limsup}{\mathop{\rm Limsup}}
\providecommand{\eLiminf}{\rm e\text{-}\mathop{\rm liminf}}
\providecommand{\eLimsup}{\rm e\text{-}\mathop{\rm limsup}}

\definecolor{forestgreenweb}{rgb}{0.13, 0.55, 0.13}

\begin{document}
	
	\begin{center}
		{\LARGE A slope generalization of Attouch theorem}
		
		\vspace{1cm}
		
		{\Large \textsc{Aris Daniilidis, David Salas, Sebasti{\'a}n Tapia-Garc{\'i}a}}
	\end{center}
	
	\bigskip
	
	\noindent\textbf{Abstract.} A classical result of variational analysis, known
	as Attouch theorem, establishes the equivalence between epigraphical
	convergence of a sequence of proper convex lower semicontinuous functions and
	graphical convergence of the corresponding subdifferential maps up to a
	normalization condition which fixes the integration constant. In this work, we
	show that in finite dimensions and under a mild boundedness assumption, we can replace subdifferentials (sets of vectors) by slopes (scalars,
	corresponding to the distance of the subdifferentials to zero) and still
	obtain the same characterization: namely, the epigraphical convergence of
	functions is equivalent to the epigraphical convergence of their slopes. This
	surprising result goes in line with recent developments on slope
	determination \cite{BCD2018, PSV2021} and slope sensitivity \cite{DD2023} for
	convex functions.
	
	\bigskip
	
	
	\vspace{0.6cm}
	
	\textit{Secondary} 37C10, 46N10, 49K40.
	
	\tableofcontents
	

	\section{Introduction}\label{sec:intro}
	
	In 1977, H\'{e}dy Attouch showed that a sequence of proper convex lower
	semicontinuous functions $\{f_{n}\}_{n\geq1}$ epi-converges to (a lower
	semicontinuous convex function) $f$ if and only if the sequence $\{\partial
	f_{n}\}_{n\geq1}$ of the corresponding subdifferentials converges graphically
	to the subdifferential $\partial f$ of $f$ and a normalization condition
	(fixing the constant of integration) holds, see~\cite{A1977, A1984-book}.
	Epi-convergence of a sequence of functions refers to the set-convergence of
	the sequence of epigraphs of the functions, while the graph convergence of a
	sequence of set-valued mappings involves the set-convergence of their graphs.
	In finite dimensions, both convergences are in the Painlev\'{e}--Kuratowski
	sense. The result remains valid in a reflexive Banach space, provided the
	convergence of $\{\mathrm{epi}f_{n}\}_{n\geq1}$ to $\mathrm{epi}f$ is taken in
	the Mosco sense (see \cite{M-1969}). \smallskip\newline Attouch theorem has
	been further extended in \cite{AB1993,CT1998} to any Banach space, using the
	notion of slice convergence which is shown to be equivalent to the Mosco
	epi-convergence of both functions and their convex conjugates. Further
	extensions cover more general classes of functions, as for instance the class
	of primal lower nice functions (see \cite{P-1992, LPT-1995, CT1998} \textit{e.g.}).
	\smallskip\newline The importance of the Attouch theorem can be measured by
	its numerous applications: it has been used to establish strong solutions in
	parabolic variational inequalities \cite{AD-1978}, stability results in
	numerical optimization \cite{L-1988} as well as theoretical results on
	generalized second order derivatives of convex functions \cite{R-1990} or in
	relation with the differentiability of Lipschitz set-valued maps
	\cite{DQ-2022}. It also meets applications in non-regular mechanics and in
	subgradient evolution systems, see \cite{A1984-book, ABM2014-book} and
	references therein.\smallskip\newline The original proofs of the Attouch
	theorem (see \cite{A1977, A1984-book, AB1993}) are based on the integration
	formula of Rockafellar \cite{R-1970} for the class of maximal cyclically
	monotone operators, which is a characteristic property of the subdifferential
	map of a convex function. The approach of \cite{CT1998} is different, but
	still relies on the \textit{subdifferential determination} of any convex
	function. Indeed, it is well-known that the equality $\partial f=\partial g$
	for any two convex lower semicontinuous functions $f,g$ guarantees that the
	functions are equal up to a constant. \smallskip\newline Quite recently, the
	following intriguing result has been established: convex lower semicontinuous
	functions are fully determined by the slope mapping $x\mapsto s_{f}(x):=\mathrm{dist}(0,\partial f(x))$ (rather than the whole subdifferential), up to an additive constant, provided they are bounded from below. In other
	words, knowledge of the remoteness of the subdifferential (which is a scalar) at every point, 
	gives in this case, enough information for the full determination of the subdifferential and consequently, 
	of the function, that is,
	\begin{equation}
		s_{f}\,=s_{g}\qquad\;\Longleftrightarrow\;\qquad\partial f=\partial
		g\qquad\;\Longleftrightarrow\;\qquad f=g+\mathrm{cst.}\label{eq:det}
	\end{equation}
	This result has first been established in Hilbert spaces for the class of smooth (convex and bounded from below) functions \cite{BCD2018} and has then been extended to the class of (nonsmooth) convex continuous and bounded from
	below functions \cite{PSV2021}. Although it is not relevant for our purposes, let us mention for completeness that \eqref{eq:det} 
	was ultimately established in~\cite{TZ2023} for convex functions defined in an arbitrary Banach space.
	Further extensions to the class of Lipschitz functions in metric spaces, upon
	knowledge of the critical set, have been done in \cite{DS2022}.\smallskip \newline 
	Very recently, a study of robustness of the \textit{slope determination} result has 
	been carried out in~\cite{DD2023}, motivated by the following question:
	
	\begin{center}
		\textit{If the slopes of two convex functions are close, are the function
			values close}\,?
	\end{center}
	
	In finite dimensions, the main result of~\cite{DD2023} reads, roughly
	speaking, as follows: if $f,g$ are two convex continuous functions that attain
	their minimum value, then:
	\[
	\Vert g-f\Vert_{\mathcal{U}}\lesssim\Vert s_{g}-s_{f}\Vert_{\mathcal{U}}+\Vert
	g-f\Vert_{C_{f}\cup C_{g}},
	\]
	where $\mathcal{U}$ is any bounded set, $\Vert\cdot\Vert_{\mathcal{U}}$ is the
	sup-norm over $\mathcal{U}$, $C_{f}:=\mathrm{argmin\,}f$ and $C_{g}
	:=\mathrm{argmin\,}g$. In particular, the quantity $\Vert g-f\Vert
	_{\mathcal{U}}$ is controlled in a Lipschitz manner by the slope deviation
	$\Vert s_{g}-s_{f}\Vert_{\mathcal{U}}$, yielding the following convergence result:
	
	\begin{theorem}
		[\cite{DD2023}, Corollary~3.3]\label{thm-dim} Let $\{f_{n}\}_{n\geq 0}$ be
		convex continuous functions such that
		\[
		\mathcal{C}_{f_{n}}:=\mathrm{argmin\,}f_{n}\neq\emptyset\quad\text{for
			all\,}\, n\geq 0 \quad\text{and}\quad\mathcal{C}:=\left(
		\cup_{n\geq 0}\,\mathcal{C}_{f_{n}}\right)  \,\text{ is bounded.}
		\]
		Assume further that:\smallskip\newline\emph{(i).} $\{s_{f_{n}}\}_{n}$
		converges to $s_{f_0}$ uniformly on bounded sets;\smallskip\newline\emph{(ii)}.
		$\{f_{n}\}_{n}$ converges to $f_0$ uniformly on $\mathcal{C}$.\medskip\newline
		Then $\{f_{n}\}_{n}$ converges to $f_0$ uniformly on bounded sets.
	\end{theorem}
	
	The assumption of existence of (global) minima in the above result is
	suboptimal, since it is stronger than mere boundedness from below, which was
	the main assumption in~\eqref{eq:det}, see also \cite[Remark~3.4]{DD2023}. In
	addition, Theorem~\ref{thm-dim} does not cover variational deviations, which
	is the natural framework of the Attouch theorem.\smallskip\newline In this
	work we generalize the result of \cite{PSV2021} (slope determination) and
	complement the result of~\cite{DD2023} (slope sensitivity), establishing a
	slope version of the Attouch theorem in finite dimensions,
	under the condition that the limiting function $f$ is bounded from below.
	Since graphic
	convergence of subdifferentials is ostensibly much stronger than
	epi-convergence of the slopes (see Section~\ref{sec: easy implication} for a
	formal proof of this implication), the converse implication is the core of our
	main result (see Section~\ref{sec: main}). Therefore, in a sensitivity
	framework, our main theorem (\textit{c.f.} Theorem~\ref{theo: main})
	generalizes the Attouch theorem, in a similar way that the slope determination
	generalizes subdifferential determination.

	\subsection{Basic setting and notation}
	
	We consider the $d$-dimensional Euclidean space $\mathbb{R}^{d}$ endowed with
	its usual inner product $\langle\cdot,\cdot\rangle$ and its Euclidean norm
	$\Vert\cdot\Vert$. For a subset $A\subset\mathbb{R}^{d}$, we denote by
	$\into(A)$, $\overline{A}$, $\partial A$ and $\ri(A)$ its interior, closure,
	boundary and relative interior, respectively. Given $x\in\mathbb{R}^{d}$,
	we write $B(x,r)$ and $\overline{B}(x,r)$ to denote the open and closed $r$-balls centered at $x$, and
	we
	define its distance to the set $A$ as follows:
	\[
	\mathrm{dist\,}(x,A):=\,\underset{a\in A}{\inf}\,\Vert x-a\Vert.
	\]
	For a function $f:\mathbb{R}^{d}\rightarrow\mathbb{R}\cup\{+\infty\}$, we
	denote its effective domain and respectively its epigraph by:
	\[
	\dom\,f:=\{x\in\mathbb{R}^{d}:f(x)<+\infty\}\quad\text{and}\quad\epi
	f=\{(x,\alpha)\in\mathbb{R}^{d}\times\mathbb{R}:\ \alpha\geq f(x)\}.
	\]
	The (Moreau-Rockafellar) subdifferential of $f$ is then defined as follows:
	\begin{equation}\label{pino}
		\partial f(x)=\{x^{\ast}\in\mathbb{R}^{d}:~f(y)\geq f(x)+\langle x^{\ast
		},y-x\rangle,~\forall y\in\mathbb{R}^{d}\},
	\end{equation}
	if $x\in\dom\,f$ and empty otherwise. Note that $f$ may not be a convex
	function and $\partial f(x)$ may be empty even if $x\in\dom\,f$. If $f$ is
	proper (i.e., $\dom\,f\neq\emptyset$), we denote by $f^{\ast}:\mathbb{R}
	^{d}\rightarrow\mathbb{R}\cup\{+\infty\}$ its Fenchel conjugate, that is,
	\begin{equation}
		f^{\ast}(x^{\ast})=\sup_{x\in\mathbb{R}^{d}}\left\{  \langle x^{\ast}
		,x\rangle-f(x)\right\}  .\label{eq:Conjugate}
	\end{equation}
	It is easy to check from the definition of $f^{\ast}$ that Young-Fenchel
	inequality holds true: for all $(x,x^{\ast})\in\mathbb{R}^{d}\times
	\mathbb{R}^{d}$ one has that $f(x)+f^{\ast}(x^{\ast})\geq\langle x^{\ast
	},x\rangle$. Moreover, the subdifferential of $f$ can be characterized in
	terms of its conjugate function as follows:
	\begin{equation}
		x^{\ast}\in\partial f(x)\iff f(x)+f^{\ast}(x^{\ast})=\langle x^{\ast}
		,x\rangle.\label{eq:CharacSubdiff-Fenchel}
	\end{equation}
	Following \cite{DgMT-1980} we define the (metric or local) slope of a function
	$f$ at a point $x\in\mathbb{R}^{d}$ as follows:
	\[
	s_{f}(x)=\left\{
	\begin{array}
		[c]{cl}
		\underset{y\rightarrow y}{\limsup}\,\frac{\{f(x)-f(y)\}^{+}}{d(y.x)}, & \text{
			if }x\in\dom\,f\\
		+\infty, & \text{ otherwise, }
	\end{array}
	\right.
	\]
	where $\{\alpha\}^{+}=\max\,\{0,\alpha\}$ and $d(y,x)=\Vert x-y\Vert$. This
	notion has been extensively studied in the framework of metric analysis, see
	\cite{I-2000, AC-2004, AGS-2008, DIL-2015, AC-2017} and references therein. In
	the special case that the function $f$ is convex and lower semicontinuous, for
	every $x\in\mathbb{R}^{d}$ one has:
	\begin{equation}
		s_{f}(x)=\mathrm{dist\,}(0,\partial f(x))\qquad\text{(distance of the
			subdifferential to }0\text{)}.\label{eq:SlopeDistanceSubdiff}
	\end{equation}
	Whenever $f$ is a proper convex lower semicontinuous function and $x\in\dom\,f$, it is
	well-known that $\partial f(x)$ is a convex closed set. Notice that
	$s_{f}(x)=+\infty$ if and only if $\partial f(x)=\emptyset.$ Whenever
	$\partial f(x)$ is nonempty, we denote by $\partial^{\circ}f(x)$ the (unique)
	element of minimal norm of $\partial f(x)$, that is,
	\begin{equation}
		\partial^{\circ}f(x)=\mathrm{proj}(0;\partial f(x)),\quad\forall
		x\in\dom\,\partial f,\label{eq:MinNorm}
	\end{equation}
	where $\mathrm{proj}(\cdot;A)$ stands for the projection to a set
	$A\subset\mathbb{R}^{d}$ and
	\[
	\dom\,\partial f=\{x\in\mathbb{R}^{d}\ :\ \partial f(x)\neq\emptyset\}
	\]
	is the effective domain of the subdifferential of $f$. 
	Notice that \eqref{pino} yields 
	\begin{equation}
		0\in\partial f(x)  \iff s_f(x)=0 \iff x \in \arg\min f \quad\text{(set of global minimizers of $f$)}
	\end{equation}
	
	As already mentioned in the introduction, the slope determines, up to a constant,
	the class of convex lower semicontinuous functions that are bounded from
	below. In a Hilbert space setting an important intermediate result, the
	so-called comparison principle, was established in~\cite{PSV2021}. This is
	recalled below.
	
	\begin{theorem}
		[Comparison principle]\label{theo: comparisonPrinciple} Let $f,g:\mathbb{R}^{d}\to\mathbb{R}\cup\{+\infty\}$ be two convex lower semicontinuous functions
		that are bounded from below. Assume that
		\begin{enumerate}
			\item[$(i)$.] $\inf f\geq\inf g$ ; and
			\item[$(ii)$.] $s_{f}(x)\geq s_{g}(x)$, for all $x\in\mathbb{R}^{d}.$
		\end{enumerate}
		Then it holds: $f\geq g$.
	\end{theorem}
	
	In what follows we identify the subdifferential $\partial f$ (which is a
	multi-valued map from $\mathbb{R}^{d}$ to $\mathbb{R}$) with its graph
	$\mathrm{gph}(\partial f)$ (which is a subset of $\mathbb{R}^{d}
	\times\mathbb{R}^{d}$) and we indistinctively switch from the notation
	$x^{\ast}\in\partial f(x)$ to the notation $(x,x^{\ast})\in\partial f$. Under
	this slight abuse of notation, we have:
	\begin{align*}
		\partial f  &  :=\{(x,x^{\ast})\in\mathbb{R}^{d}\times\mathbb{R}
		^{d}\,:~x^{\ast}\in\partial f(x)\}\,\,\big(\subset\mathbb{R}^{d}
		\times\mathbb{R}^{d}\,\big)\,,\\
		\triangle f  &  :=\{(x,x^{\ast},\alpha)\in\mathbb{R}^{d}\times\mathbb{R}^{d}
		\times\mathbb{R}:~x^{\ast}\in\partial f(x),~\alpha=f(x)\}\,\,\big(\subset\,
		\mathbb{R}^{d}\times\mathbb{R}^{d}\times\mathbb{R}\,\big)\,.
	\end{align*}
	
	For a proper convex lower semicontinuous function $f:\R^{d}\to \R\cup\{+\infty\}$ and a point $x_0\in \overline{\mathrm{\dom}}\, f$, we say that an absolutely continuous curve $\gamma:[0,+\infty)\to\R^{d}$ is a \textit{(maximal) steepest descent curve} for $f$ emanating from $x_0$ if it solves the differential inclusion
	\begin{equation}\label{eq:GradientFlow}
		\begin{cases}
			\, \dot{\gamma}(t) \in -\partial f(\gamma(t)),\quad\forall t\in[0,+\infty),\\
			\, \gamma(0) = x_0.
		\end{cases} 
	\end{equation}
	It is well known  (see, e.g., \cite[Chapter 17]{ABM2014-book}) that  for any initial point $x_0\in \overline{\mathrm{\dom}}\, f$, there exists a unique steepest descent curve emanating from $x_0$. In addition, the functions $t\mapsto f(\gamma(t))$ and $t\mapsto s_f(\gamma(t))$ are decreasing and satisfy 
	\begin{equation}
		\lim_{t\to+\infty} f(\gamma(t)) = \inf f \qquad \text {and } \qquad \lim_{t\to+\infty} s_f(\gamma(t)) = 0.
	\end{equation}
	If $\gamma$ is bounded (that is, $\gamma([0,+\infty))\subset B(0,M)$ for some $M>0$), then it has finite length (see \cite{MP1991, DDD2015} \textit{e.g.}). This happens exactly when $\arg\min f\neq\emptyset$ and in this case, $\gamma(t) \underset{t\to+\infty}{\longrightarrow}\gamma_{\infty}$, with $s_f(\gamma_{\infty})=0$. (Notice here that it is possible to have convergence in finite time, \textit{i.e.} $\gamma(T)=\gamma_{\infty}$ for some $T>0$, case in which $\gamma$ becomes stationary afterwards: think for example of the function $f(x)=\|x\|$, for all $x\in \mathbb{R}$.)
	
	
	\subsection{Notions of convergence and Attouch theorem}\label{subsec:NotionsOfConvergence}
	
	\label{subsec:intro} Let $\{S_{n}\}_{n}$ be a sequence of subsets of
	$\mathbb{R}^{d}$. We consider the inferior and superior limits of
	$\{S_{n}\}_{n}$ in the sense of Painlev\'{e}-Kuratowski as
	\begin{align*}
		\underset{n\rightarrow\infty}{\Liminf\,}S_{n} &  :=\left\{  x\in\mathbb{R}
		^{d}\ :\ \limsup_{n\rightarrow\infty}\mathrm{dist\,}(x,S_{n})=0\right\}  ,\\
		\underset{n\rightarrow\infty}{\Limsup\,}S_{n} &  :=\left\{  x\in\mathbb{R}
		^{d}\ :\ \liminf_{n\rightarrow\infty}\mathrm{dist\,}(x,S_{n})=0\right\}  .
	\end{align*}
	We say that $\{S_{n}\}_{n}$ converges to a set $S$ in the sense of
	Painlev\'{e}-Kuratowski, which we denote by $S_{n}\xrightarrow{PK}S$, if both,
	the inferior and superior limits of $\{S_{n}\}_{n}$ coincide with $S$. Noting
	that $\Liminf S_{n}\subset\Limsup S_{n}$, one can write
	\begin{equation}
		S_{n}\xrightarrow{PK}S\iff\Limsup_{n\rightarrow\infty}S_{n}\subset
		S\subset\Liminf_{n\rightarrow\infty}S_{n}.\label{eq:characPKConvergence}
	\end{equation}
	In what follows, given a sequence of functions $\{\phi_{n}\}_{n}$ from
	$\mathbb{R}^{d}$ to $\mathbb{R}\cup\{+\infty\}$, we define the functions
	$\phi_{l},\phi_{u}:\mathbb{R}^{d}\rightarrow\mathbb{R}\cup\{\pm\infty\}$ as the
	lower and, respectively, the upper epigraphical limits of $\{\phi_{n}\}_{n}
	$, given as follows:
	\begin{align}
		\phi_{l}(x) &  =(\underset{n\rightarrow\infty}{{\eLiminf}}\,\phi_{n}
		)(x):=\inf_{x_{n}\rightarrow x}\liminf_{n\rightarrow\infty}\,\phi_{n}
		(x_{n}),\label{eq:epi-lim}\\
		\phi_{u}(x) &  =(\underset{n\rightarrow\infty}{{\eLimsup}}\,\phi_{n}
		)(x):=\inf_{x_{n}\rightarrow x}\limsup_{n\rightarrow\infty}\,\phi_{n}
		(x_{n})\nonumber
	\end{align}
	where, in both cases, the infimum is taken over all sequences $\{x_{n}\}_{n}\subset\mathbb{R}^{d}$ 
	converging to $x$. Given a strictly increasing sequence of natural numbers $\{k(n)\}_{n}$ (which we 
	indistinctively also denote by $\{k_{n}\}_{n}$) we denote by $\phi_{l,k(n)}$ (respectively, $\phi_{u,k(n)}$) 
	the lower (respectively, upper) epigraphical limit of the subsequence $\{\phi_{k(n)}\}_{n}$ of $\{\phi_{n}\}_{n}$.
	
	\begin{remark}
		[attainability of infimum and lower semicontinuity]
		\label{rem: infsAttained} The infima that define $\phi_{l}$ and $\phi_{u}$ in~\eqref{eq:epi-lim} are attained, that is, for every $x\in\mathbb{R}^{d}$ there
		exist (infimizing) sequences $\{x_{n}^{l}\}_{n}$ and $\{x_{n}^{u}\}_{n}$,
		converging to $x$, satisfying:
		\[
		\phi_{u}(x)=\limsup_{n\rightarrow\infty}\,\phi_n(x_{n}^{u})\quad\text{ and
		}\quad\phi_{l}(x)=\liminf_{n\rightarrow\infty}\,\phi_{n}(x_{n}^{l}).
		\]
		Based on the above remark and using a diagonal argument, we easily deduce that
		the functions $\phi_{l}$ and $\phi_{u}$ are lower semicontinuous.
	\end{remark}
	
	Finally, we say that a sequence of functions $\{\phi_{n}\}_{n}$
	\emph{converges epigraphically} to a function $\phi$ and denote $\phi
	_{n}\xrightarrow{e}\phi$, if the sequence of epigraphs $\{\mathrm{epi\,}
	\phi_{n}\}_{n}$ converges to $\mathrm{epi}$\thinspace$\phi$ in the sense of
	Painlev\'{e}-Kuratowski, that is:\newline
	\[
	\phi_{n}\xrightarrow{e}\phi\qquad\iff\qquad\epi\phi_{n}
	\xrightarrow{PK}\epi\phi.
	\]
	It is well-known that
	\begin{equation}
		\phi_{n}\xrightarrow{e}\phi\iff\phi_{u}=\phi=\phi_{l}\iff \phi_{l}\geq\phi \geq\phi_{u}.\label{eq:characEpiConvergence}
	\end{equation}
	Let us finally recall, in the finite dimensional setting, the following
	celebrated variational approximation result due to H. Attouch \cite{A1977}:
	
	\begin{theorem}
		[Attouch theorem]\label{theo: attouch}Let $f,~\{f_{n}\}_{n}:\mathbb{R}^{d}\rightarrow\mathbb{R}\cup\{+\infty\}$ 
		be proper convex lower semicontinuous functions. The following assertions are equivalent:
		
		\begin{enumerate}
			\item[$(i)$.] $\epi f_{n}\xrightarrow{PK}\epi f$ \ (that is, $f_{n}
			\xrightarrow{e}f$).
			
			\item[$(ii)$.] $\partial f_{n}\xrightarrow{PK}\partial f$ and:
			\begin{equation}
				\exists(x,x^{\ast})\in\partial f\text{ and a sequence }(x_{n},x_{n}^{\ast}
				)\in\partial f_{n},\quad(x_{n},x_{n}^{\ast},f_{n}(x_{n}))\rightarrow
				(x,x^{\ast},f(x)).\tag{NC}\label{eq:NC}
			\end{equation}

			\item[$(iii)$.] $\triangle f_{n}\xrightarrow{PK}\triangle f$.
		\end{enumerate}
	\end{theorem}
	
	The normalization condition \eqref{eq:NC} is necessary in order to fix a
	reference point. Without this condition, simple counterexamples can be
	constructed: indeed, consider the functions $f_{n}(x)\equiv n,$ for all
	$n\geq1$ and the function $f(x)\equiv0.$ Then $\partial f_{n}(x)=\partial
	f(x)=\{0\},$ for all $x\in\mathbb{R}^{d}$ and $n\geq1,$ but $f_{n}
	(x)\rightarrow\infty,$ for all $x\in\mathbb{R}^{d}.$\smallskip\newline Our
	objective in this work is to provide a version of the Attouch theorem which is
	based on the epigraphical convergence of the slope mappings $s_{f_{n}
	}:\mathbb{R}^{d}\rightarrow\mathbb{R}\cup\{+\infty\}$ (rather than the
	graphical convergence of the subdifferential maps $\partial f_{n}
	:\mathbb{R}^{d}\rightrightarrows\mathbb{R}^{d}$). Before we proceed to this,
	let us extract the following consequence of Theorem~\ref{theo: attouch} (i)$\Rightarrow$(ii) for future use.
	\medskip
	
	\begin{remark}
		\label{rem_key-rem}Let $f_{n}\xrightarrow{e}f.$ Then for every strictly increasing sequence $\{k_n\}_{n\geq 1}$
		and for every $\{(x_{k_n}, x_{k_n}^{\ast})\}_{n}\subset\mathbb{R}^{d}\times\mathbb{R}^{d}$ such that
		$x_{k_n}\rightarrow x,$ $x_{k_n}^{\ast}\rightarrow x^{\ast}$ and $x_{k_n}^{\ast}
		\in\partial f_{k_n}(x_{k_n}),$ we have $x^{\ast}\in\partial f(x).$
	\end{remark}

	\subsection{Our contribution}
	
	The goal of this work is to establish that epigraphical convergence of convex
	functions can be characterized by epigraphical convergence of the slopes. \smallskip\newline
	Our approach relies on the determination result of \cite{PSV2021} and naturally
	inherites the restriction that the limit function should be bounded from
	below. As in the Attouch theorem, a normalization condition will also be
	required. In this work we can either use the same condition \eqref{eq:NC} as
	in Theorem~\ref{theo: attouch} or an alternative condition over the infimum of
	the epigraphical lower and upper limits. Concretely, our main result is as follows:
	
	\begin{theorem}
		[main result]\label{theo: main} Let $f,\{f_{n}\}_{n}:\mathbb{R}^{d}
		\rightarrow\mathbb{R}\cup\{+\infty\}$ be proper convex lower semicontinuous
		functions. Assume that $\inf f\in\mathbb{R}$. Then, the following assertions
		are equivalent:
		
		\begin{enumerate}
			\item[(i).] $f_{n}\xrightarrow{e}f$.
			
			\item[(ii).] $s_{f_{n}}\xrightarrow{e}s_{f}$ and \eqref{eq:NC} holds.
			
			\item[(iii).] $s_{f_{n}}\xrightarrow{e}s_{f}$ and $\inf f_{l}=\inf f=\inf
			f_{u}$.
		\end{enumerate}
	\end{theorem}
	
	The rest of the manuscript is organized as follows: in
	Section~\ref{sec: easy implication} we show that implications~(ii) and~(iii)
	of the statement of Theorem~\ref{theo: main} follow easily from~(i) and
	Theorem~\ref{theo: attouch} (Attouch theorem). Then,
	Section~\ref{sec: useful results}, is devoted to a preliminary study of the
	functions $f_{l}$ and $f_{u}$. Finally, in Section~\ref{sec: main}, we show
	that either one of (ii) or (iii) implies (i). The approach is divided into two
	parts: we first show that $f_{u}\leq f$ in Subsection~\ref{subsec: fu leq f},
	and then in Subsection~\ref{subsec: gap} we prove that $f\leq f_{l}$. The main
	result and final comments are given at the end (Section~\ref{subsec: Ortega}).
	

	\section{From epigraphical convergence to slope convergence}
	
	\label{sec: easy implication}
	
	In this section, we show the ``easy"
	\ implications of Theorem~\ref{theo: main}, namely, (i)$\Rightarrow
	$(ii),(iii). The proof consists of studying the upper and the lower
	epigraphical limits of the slope sequence $\{s_{f_{n}}\}_{n}$ then combine
	with \eqref{eq:characEpiConvergence} to deduce the result. A standard argument
	that will repeatedly appear in this work, is to study separately the points
	where the limit function (in this case $s_{f}$) is finite from those where the
	limit is infinite.
	
	\begin{theorem}
		\label{theo: easy implication} Let $f,\{f_{n}\}_{n}:\mathbb{R}^{d}
		\rightarrow\mathbb{R}\cup\{+\infty\}$ be proper convex lower semicontinuous
		functions. Assume that $f_{n}\xrightarrow{e}f$. Then
		\[
		s_{f_{n}}\xrightarrow{e}s_{f},\qquad\inf\,f=\inf\,f_{l}=\inf\,f_{u}
		\qquad\text{and}\qquad\eqref{eq:NC}\text{ holds.}
		\]
		
	\end{theorem}
	
	\noindent\textit{Proof.} Our assumption yields $f=f_{u}=f_{l}$, thus, $\inf\,f=\inf\,f_{l}=\inf\,f_{u}$. Condition~$\eqref{eq:NC}$
	follows from Theorem~\ref{theo: attouch} (i)$\Rightarrow$(ii). It remains
	to prove the epiconvergence of the sequence $\{s_{f_{n}}\}_{n}$ to $s_{f}$. To
	this end, let $x\in\dom\,f$ and consider separately two cases:
	
	\begin{itemize}
		\item \textit{Case 1}: $\quad\partial f(x)=\emptyset$ (that is, $s_{f}(x)=+\infty$)
	\end{itemize}
	
	In this case, we need to show that $({\eLiminf}\,s_{f_{n}})(x)=+\infty$. Let
	us assume, towards a contradiction, that there exists a sequence
	$\{x_{n}\}_{n}\subset\mathbb{R}^{d}$ converging to $x$, such that
	${\underset{n\rightarrow\infty}{\lim\inf\,}s_{f_{n}}(x_{n})<+\infty.}$ Then,
	for an adequate subsequence $\{x_{k_{n}}\}_{n}$ we would have
	\[
	\underset{n\rightarrow\infty}{\lim\inf}\,s_{f_{n}}(x_{n})=\,\underset{n\rightarrow\infty}{\lim}\,s_{f_{k_{n}}}(x_{k_{n}})<\infty\,,
	\]
	and (up to a new subsequence) $x_{k_{n}}^{\ast} \rightarrow x^{\ast}$, for some $ x^{\ast}\in \mathbb{R}^d$, where 
	$x_{k_{n}}^{\ast}:=\partial^{\circ}f_{k_{n}}(x_{k_{n}})$ is the element of minimal norm in $\partial f_{k_{n}}(x_{k_{n}})$, as in~\eqref{eq:MinNorm}.
	By Remark~\ref{rem_key-rem} we infer that $x^{\ast}\in\partial f(x)$, which is a contradiction. Therefore, $ ({\eLiminf}\,s_{f_{n}})(x)=+\infty=s_{f}(x)$.
		
	\begin{itemize}
		\item \textit{Case 2}: $\quad\partial f(x)\neq\emptyset$ (that is,
		$s_{f}(x)<+\infty$)
	\end{itemize}
	
	Let $x\in\dom s_{f}$ and $\bar{x}^{\ast}\in\partial f(x)$ such that $\Vert
	\bar{x}^{\ast}\Vert=s_{f}(x)$. Since $(x,\bar{x}^{\ast},f(x))\in\triangle f$
	and since $\triangle f_{n}\xrightarrow{PK}\triangle f$ (\textit{c.f.}
	Theorem~\ref{theo: attouch}), there exists a sequence ${(x_{n},x_{n}^{\ast
		},f_{n}(x_{n}))\in\triangle f_{n}}$ converging to $(x,\bar{x}^{\ast},f(x))$.
	Thus, using \eqref{eq:SlopeDistanceSubdiff}, we deduce
	\[
	\limsup_{n\rightarrow\infty}\,s_{f_{n}}(x_{n})\,\leq\,\limsup_{n\rightarrow
		\infty}\,\Vert x_{n}^{\ast}\Vert=\Vert\bar{x}^{\ast}\Vert=s_{f}(x),
	\]
	which yields $({\eLimsup}\,s_{f_{n}})(x)\leq s_{f}(x)$. It remains to show that
	\[
	\inf_{x_{n}\rightarrow x} \liminf_{n\rightarrow\infty}\,s_{f_{n}}(x_{n})\geq
	s_{f}(x).
	\]
	To this end, we consider an arbitrary sequence $\{x_{n}\}_{n}\subset
	\mathbb{R}^{d}$ converging to $x$. For a suitable subsequence $\{k_{n}
	\}_{n}\ $we have:
	\[
	\lim_{n\rightarrow\infty}\,s_{f_{k_{n}}}(x_{k_{n}})\,=\,\liminf_{n\rightarrow
		\infty}\,s_{f_{n}}(x_{n})\,=\rho
	\]
	and we need to show that $\rho\geq s_{f}(x).$ We can obviously assume that
	$\rho<+\infty$. For each $n\in\mathbb{N}$, let $x_{k_{n}}^{\ast}=\partial
	^{\circ}f_{k_{n}}(x_{k_{n}})$ be the element of minimal norm of $\partial
	f_{k_{n}}(x_{k_{n}})$, that is, $\Vert x_{k_{n}}^{\ast}\Vert=s_{f_{k_{n}}
	}(x_{k_{n}})$. Since $s_{f_{k_{n}}}(x_{k_{n}})\rightarrow\rho,$ the
	subsequence $\{x_{k_{n}}^{\ast}\}_{n}$ is bounded and converges (up to a new
	subsequence) to some $x^{\ast}\in\mathbb{R}^{d}$, with $\Vert x^{\ast}
	\Vert=\rho$. By Remark~\ref{rem_key-rem} we have $x^{\ast}\in\partial f(x)$
	and consequently
	\[
	s_{f}(x)\leq\Vert x^{\ast}\Vert\,=\,\lim_{n\rightarrow\infty}\,\Vert x_{k_{n}
	}^{\ast}\Vert\,=\lim_{n\rightarrow\infty}\,s_{f_{k_{n}}}(x_{k_{n}})=\rho.
	\]
	Since the sequence $\{x_{n}\}_{n}$ is arbitrary, we have $s_{f}(x)\leq
	({\eLiminf}\,s_{f_{n}})(x)$.\medskip
	
	The proof is complete. \hfill$\Box$
	

	
	\section{Some intermediate results}
	
	\label{sec: useful results}
	

	In this part, we obtain some preliminary results, which are needed for the
	proof of the \textquotedblleft difficult\textquotedblright\ implication of our
	main theorem. Some of the forthcoming results are essentially known, other are
	less obvious and require a careful analysis.
	
	\subsection{General results from convex analysis}
	
	The first result is essentially known.
	
	\begin{proposition}
		\label{prop: convex rint} Let $f:\mathbb{R}^{d}\rightarrow\mathbb{R}
		\cup\{+\infty\}$ be a proper convex lower semicontinuous function. Then
		$\ri(\dom s_{f})=\ri(\dom f)$ and therefore it is a convex set.
	\end{proposition}
	
	\noindent\textit{Proof.} Let us first notice that $\dom s_{f}=\dom\partial f.$
	Since $f$ is convex, $\partial f$ is nonempty on $\ri(\dom\,f)$. Thus,
	$\ri(\dom\,f)\subset\dom s_{f}\subset\dom\,f$. Without loss of generality, we
	may assume that $0\in\dom\,f$. Set $V=\mathrm{span}(\dom\,f)$, that is, the
	subspace of $\mathbb{R}^{d}$ generated by $\dom\,f$. Notice that
	$\ri(\dom\,f)$ generates the same subspace $V$, therefore $\mathrm{span}
	(\dom
	s_{f})=V$. Since the relative interiors of $\dom s_{f}$ and $\dom\,f$ are
	taken with respect to the same space $V$, we have $\ri(\dom s_{f}
	)\subset\ri(\dom\,f)$. The conclusion follows from the convexity of
	$\dom\,f$.$\,$\hfill$\Box$
	
	\bigskip
	
	The following result is also quite intuitive.
	
	\begin{proposition}
		\label{prop: extension to the adherence} Let $f,g:\mathbb{R}^{d}\rightarrow\mathbb{R}\cup\{+\infty\}$ be two functions, with $f$ convex and $g$ lower semicontinuous. Let $A\subset\dom\,f$ be a nonempty convex set. Assume that $f\geq
		g$ on $A$. Then, $f\geq g$ on $\overline{A}$.
	\end{proposition}
	
	\noindent\textit{Proof.} Let $\overline{x}\in\overline{A}$ and $x\in\ri(A).$
	Then $(\overline{x},x]\subset A$. Note that $f(x)\geq g(x)\in\mathbb{R}$.
	Since $f$ is convex and $g$ is lower semicontinuous, we have that
	\begin{align*}
		g(\overline{x}) \leq\liminf_{t\rightarrow0^{+}}\,g(tx+(1-t)\overline{x}
		)\leq\liminf_{t\rightarrow0^{+}}\,f(tx+(1-t)\overline{x}) \leq\liminf
		_{t\rightarrow0^{+}}\,\left\{  tf(x)+(1-t)f(\overline{x})\right\}
		=f(\overline{x}).
	\end{align*}
	Since $\overline{x}$ was arbitrarily chosen in $\overline{A},$ we conclude
	that $f\geq g$ on $\overline{A}$.\hfill$\Box$
	
	\bigskip
	
	Let $K\subset\mathbb{R}^{d}$ be a nonempty convex set. We denote by $\sigma
	_{K}:\mathbb{R}^{d}\rightarrow\mathbb{R}$ the support function of $K$, that
	is, for any $x\in\mathbb{R}^{d}$ we have
	\[
	\sigma_{K}(x):=\,\sup_{y\in K}\,\langle x, y\rangle.
	\]
    Additionally, for $x\in K$, we denote by $N_K(x)$ the normal cone of $K$ at $x$. It is well known that
	\[
	N_K(x) = \{ x^*\in \R^{d}\, :\,  \sigma_{K}(x^{\ast})\leq \langle x^{*},x\rangle \}.
	\]
With this in mind, the following proposition establishes a density characterization for the subdifferential of convex functions.
	\begin{proposition}
		\label{Prop: subdiff r-int} Let $f:\mathbb{R}^{d}\rightarrow\mathbb{R}
		\cup\{+\infty\}$ be a proper lower semicontinous convex function and let
		$(x,x^{\ast})\in\dom\,f\times\mathbb{R}^{d}$. Assume there exists a dense
		subset $D$ of $\ri(\dom\,f)$ such that
		\begin{equation}
			\forall y\in D,\,\exists y^{\ast}\in\partial f(y)\text{ such that }\langle
			y^{\ast}-x^{\ast},y-x\rangle\geq0.\label{eq: ConditionOnRI-subdiff}
		\end{equation}
		Then $x^{\ast}\in\partial f(x)$.
	\end{proposition}
	
	\noindent\textit{Proof.} Without loss of generality, we may assume
	$0\in\ri(\dom\,f)$ and set $V=\mathrm{span}(\dom\,f)$. We consider two cases.
	
	\begin{itemize}
		\item \textit{Case 1}: $\;V=\mathbb{R}^{d}$ and consequently
		$\ri(\dom\,f)=\mathrm{int}(\dom\,f)$
	\end{itemize}
	
	In this case, $\partial f$ is locally bounded on $\mathrm{int}(\dom\,f)$ and
	upper semicontinuous (in the sense of set-valued maps). Therefore
	\eqref{eq: ConditionOnRI-subdiff} entails that
	\[
	\forall y\in\into(\dom\,f),\,\exists y^{\ast}\in\partial f(y)\text{ such that
	}\langle y^{\ast}-x^{\ast},y-x\rangle\geq0.
	\]
	In particular, for every differentiability point $y$ of $f$, it holds that
	$\langle\nabla f(y)-x^{\ast},y-x\rangle\geq0$. Moreover, since for every
	$y\in\into(\dom\,f)$ we have that $\partial f(y)$ can be recovered as the
	convex hull of limits of gradients of $f$ (see, e.g., \cite[Theorem
	9.61]{RW1998}), we get that
	\[
	\langle y^{\ast}-x^{\ast},y-x\rangle\geq0,\quad\forall(y,y^{\ast}
	)\in\mathcal{D},
	\]
	where
	\[
	\mathcal{D}=\overline{\{(z,z^{\ast}):\ z\in\into(\dom\,f),\,z^{\ast}
		\in\partial f(z)\}}.
	\]
	Now, take $(\bar{y},\bar{y}^{\ast})\in\partial f\setminus\mathcal{D}$.
	Clearly, $\bar{y}\in\partial(\dom\,f)$ and $f(\bar{y})\in\mathbb{R}$. Consider
	the convex body $K=\overline{\dom}f$ and set
	\[
	K_{n}=\left(1-\frac{1}{n}\right)K\quad\text{ and }\quad f_{n}=f+I_{K_{n}},\quad\forall
	n\in\mathbb{N},
	\]
	where $I_{K_{n}}$ is the indicator function of $K_{n},$ that is, $I_{K_{n}
	}(x)=0$, if $x\in K_{n}$ and $+\infty$ elsewhere. Notice that $K_{n}
	\subset\ri(K)$ and that by construction $f_{n}\xrightarrow{e}f$. Applying
	Theorem~\ref{theo: attouch} (Attouch theorem), we deduce that $\Delta
	f_{n}\xrightarrow{PK}\Delta f$. Therefore, there exists a sequence
	$\{(z_{n},z_{n}^{\ast},f_n(z_{n}))\}_{n}$ converging to $(\bar{y},\bar{y}^{\ast
	},f(\bar{y}))$ such that $z_{n}\in K_{n}$ and $z_{n}^{\ast}\in\partial
	f_{n}(z_{n})$. By the sum rule for subdifferentials (see, e.g., \cite[Theorem~3.16]{Phelps1993}), there exist $y_{n}^{\ast}\in\partial f(z_{n})$ and
	$v_{n}^{\ast}\in\partial I_{K_{n}}(z_{n})\equiv N_{K_{n}}(z_{n})$ such that
	$z_{n}^{\ast}= y_{n}^{\ast}+v_{n}^{\ast}.$ In particular, $(z_n,y_n^{\ast})\in\mathcal{D}$ and therefore $\langle y_n^{\ast}-x^{\ast},z_n-x\rangle\geq 0$.\smallskip\newline 
	
	Notice further
	that $f_{n}(z_{n})=f(z_{n})\rightarrow f(\bar{y})$ and that $\sigma_{K_{n}
	}=(1-\tfrac{1}{n})\sigma_{K}$. Therefore, for every $z\in K_{n}$ we have that
	\begin{align*}
		N_{K_{n}}(z)  &  =\left\{  v^{\ast}\in\mathbb{R}^{d}\ :\ \sigma_{K_{n}
		}(v^{\ast})\leq\langle v^{\ast},z\rangle\right\} \\
		&  =\left\{  v^{\ast}\in\mathbb{R}^{d}\ :\ \sigma_{K}(v^{\ast})\leq\langle
		v^{\ast},(1-\tfrac{1}{n})^{-1}z\rangle\right\}  =N_{K}\left((1-\tfrac{1}{n})^{-1}z\right).
	\end{align*}
	In particular, since $v_{n}^{\ast}\in N_{K_{n}}(z_{n})=N_{K}((\tfrac{n}
	{n-1})z_{n}),$ for every $x\in K$ we have
	\[
	\langle v_{n}^{\ast},x-(\tfrac{n}{n-1})z_{n}\rangle\leq0
	\]
	and combining with the definition of the subdifferential map given in \eqref{pino}, we deduce that
	\begin{align*}
		\langle z_{n}^{\ast}-x^{\ast},z_{n}-x\rangle &  =\langle y_{n}^{\ast}-x^{\ast
		},z_{n}-x\rangle+\langle v_{n}^{\ast},z_{n}-x\rangle\geq\langle v_{n}^{\ast
		},z_{n}-x\rangle\\
		&  =\langle v_{n}^{\ast},(\tfrac{n}{n-1})z_{n}-x\rangle-\frac{1}{n-1}\langle
		v_{n}^{\ast},z_{n}\rangle\geq-\frac{1}{n-1}\langle v_{n}^{\ast},z_{n}\rangle\\
		&  =\frac{1}{n-1}(\langle y_{n}^{\ast},z_{n}-0\rangle-\langle z_{n}^{\ast
		},z_{n}\rangle) \geq\frac{1}{n-1}\left(  f(z_{n})-f(0)-\langle z_{n}^{\ast
		},z_{n} \rangle\right)  .
	\end{align*}
	Since $(z_n,z^{\ast}_n,f(z_n)) \rightarrow (\bar y, \bar y^{\ast}, f(\bar y))$, we 
	can take limit at both sides of the obtained inequality, deducing that
	$\langle\bar{y}^{\ast}-x^{\ast},\bar{y}-x\rangle\geq0$. Since $(\bar{y}
	,\bar{y}^{\ast})$ was arbitrarily chosen in $\partial f,$ we obtain:
	\[
	\langle y^{\ast}-x^{\ast},y-x\rangle\geq0,\quad\forall(y,y^{\ast})\in\partial
	f,
	\]
	and we conclude that $(x,x^{\ast})\in\partial f$ by maximal monotonicity of
	the subdifferential (see, e.g., \cite[Theorem 3.24]{Phelps1993}).
	
	\begin{itemize}
		\item \textit{Case 2}: \ (general case)
	\end{itemize}
	
	Let $\pi_{V}:\mathbb{R}^{d}\rightarrow V$ be the orthogonal projection
	onto $V$ and let us denote by $g:V\rightarrow\mathbb{R}\cup\{+\infty\}$ the restriction
	of $f$ on $V$. Note that for every $z\in\dom g$, every $z^{\ast}\in V$ and every $\nu^{\ast}\in V^{\perp}$, we have:
	\begin{align*}
		z^{\ast}\in\partial g(z) &  \iff\forall z^{\prime}\in V,\;\langle z^{\ast},z^{\prime}-z\rangle\leq g(z^{\prime})-g(z)\\
		&  \iff\forall z^{\prime}\in V,\;\langle z^{\ast}+\nu^{\ast},z^{\prime}
		-z\rangle\leq f(z^{\prime})-f(z)\\
		&  \iff\forall z^{\prime}\in\mathbb{R}^{d},\;\langle z^{\ast}+\nu^{\ast
		},z^{\prime}-z\rangle\leq f(z^{\prime})-f(z)\;\iff\;z^{\ast}+\nu^{\ast}
		\in\partial f(z),
	\end{align*}
	that is,
	\[
	\partial f(z)=\pi_{V}^{-1}(\partial g(z)),\quad\forall z\in\dom\,f.
	\]
	Therefore, it is enough to verify that $\pi_{V}(x^{\ast})\in\partial g(x)$.
	Note that, by projecting the subgradients of $f$ onto $V$,
	\eqref{eq: ConditionOnRI-subdiff} entails that
	\[
	\forall y\in D,\exists y^{\ast}\in\partial g(y)\text{ such that }\langle
	y^{\ast}-\pi_{V}(x^{\ast}),y-x\rangle\geq0,
	\]
	where $D$ is a dense subset of $\into(\dom g)$ and where the interior is taken
	with respect to the space~$V$. Using the same reasoning as in Case 1 above, we
	conclude that $\pi_{V}(x^{\ast})\in\partial g(x)$. The proof is complete.
	\hfill$\Box$
	

	\subsection{Two key lemmas}
	
	We now state and prove two important technical lemmas. The first one
	states that for a sequence of proper convex lower semicontinuous functions
	$\{f_{n}\}_{n}$ from $\mathbb{R}^{d}$ to $\mathbb{R}\cup\{+\infty\}$, if a sequence of points~$\{x_{n}\}_{n}$ converges and the sequence of slopes $\{s_{f_{n}}(x_{n})\}_{n}$
	is bounded, then the sequence~$\{x_{n}\}_{n}$ automatically infimizes the
	expressions of $f_{u}$ and $f_{l}$ given in (\ref{eq:epi-lim}) evaluated at
	its point of convergence.
	
	\begin{lemma}
		\label{lemma: everySeqTight} Let $f_{n}:\mathbb{R}^{d}\rightarrow
		\mathbb{R}\cup\{+\infty\},$ $n\geq1,$ be a sequence of proper lower
		semicontinuous convex functions. Let $\{x_{n}\}_{n}\subset\mathbb{R}^{d}$ be such
		that $\{s_{f_{n}}(x_{n})\}_{n}$ is bounded. Assume that $\{x_{n}\}_{n}$
		converges to some $\bar{x}$. Then
		\[
		f_{l}(\bar{x})=\liminf_{n\rightarrow\infty}\,f_{n}(x_{n})\qquad\text{and}
		\qquad f_{u}(\bar{x})=\limsup_{n\rightarrow\infty}\,f_{n}(x_{n}).
		\]
		
	\end{lemma}
	
	\noindent\textit{Proof.} Let $\{y_{n}\}_{n}\subset\mathbb{R}^{d}$ be an
	arbitrary sequence converging to $\bar{x}$. Then
	\[
	\liminf_{n\rightarrow\infty}\,f_{n}(y_{n})\,\geq\,\liminf_{n\rightarrow\infty
	}\,\left\{  \,f_{n}(x_{n})-s_{f_{n}}(x_{n})\,\Vert x_{n}-y_{n}\Vert\,\right\}
	\,=\,\liminf_{n\rightarrow\infty}\,f_{n}(x_{n}).
	\]
	It follows readily that $\underset{n\rightarrow\infty}{\lim\inf}\,f_{n}
	(x_{n})=f_{l}(\bar{x})$.\smallskip\newline Let further $\{k_{n}\}_{n}$ be a
	strictly increasing sequence such that
	\[
	\limsup_{n\rightarrow\infty}\,f_{n}(x_{n})=\lim_{n\rightarrow\infty}
	\,f_{k_{n}}(x_{k_{n}}).
	\]
	Then for every sequence $\{y_{n}\}_{n}$ converging to $\bar{x}$ we have:
	\begin{align*}
		\limsup_{n\rightarrow\infty}\,f_{n}(y_{n})\,  &  \geq\,\limsup_{n\rightarrow
			\infty}\,f_{k_{n}}(y_{k_{n}})\,\geq\,\limsup_{n\rightarrow\infty}\,\left\{
		f_{k_{n}}(x_{k_{n}})-s_{f_{k_{n}}}(x_{k_{n}})\,\Vert x_{k_{n}}-y_{k_{n}}
		\Vert\right\} \\
		&  =\,\lim_{n\rightarrow\infty}\,f_{k_{n}}(x_{k_{n}})\,=\,\limsup
		_{n\rightarrow\infty}\,f_{n}(x_{n}).
	\end{align*}
	Therefore we conclude that $\underset{n\rightarrow\infty}{\limsup\,}
	f_{n}(x_{n})=f_{u}(\bar{x})$ and the proof is complete.\hfill$\Box$
	
	\bigskip
	
	The previous result will be now used to establish our second important
	technical lemma:
	
	\begin{lemma}
		\label{lemma: domInclusions} Let $f,\{f_{n}\}_{n}:\mathbb{R}^{d}
		\rightarrow\mathbb{R}\cup\{+\infty\}$ be proper convex lower semicontinuous
		functions such that the sequence of slope functions $\{s_{f_{n}}\}_{n\geq1}$
		epigraphically converges to $s_{f}$. Assume further that there exists a
		sequence $\{x_{n}\}_{n}\subset\mathbb{R}^{d}$, $\bar{x}\in\mathrm{dom\,}
		s_{f}$ and $\alpha\in\R$ such that
		\begin{equation}\label{eq:Ortega_5000}
		\{s_{f_{n}}(x_{n})\}_{n}\,\text{\ is bounded}\qquad\text{and}\qquad
		\lim_{n\rightarrow\infty}\,\left(  x_{n},f_{n}(x_{n})\right)  =(\bar{x}
		,\alpha).
		\end{equation}
		Then,
		\[
		\mathrm{dom\,}s_{f}\,\subset\,\mathrm{dom\,}f_{l}\cap\mathrm{dom\,}
		f_{u}\,\subset\,\mathrm{dom\,}f_{l}\cup\mathrm{dom\,}f_{u}\,\subset
		\,\overline{\mathrm{dom}}\,s_{f}=\overline{\mathrm{dom}}\,f.
		\]
		
	\end{lemma}
	
	\noindent\textit{Proof.} The equality $\overline{\mathrm{dom}}\,s_{f}=\overline{\mathrm{dom}}\,f$
	follows easily from the fact that $\mathrm{dom\,}s_{f}=\mathrm{dom\,}\partial f$ is dense in $\mathrm{dom\,}f$ 
	(see, e.g., \cite[Theorem~3.17]{Phelps1993}). In addition, thanks to
	Lemma~\ref{lemma: everySeqTight}, we have $f_{l}(\bar{x})=f_{u}(\bar
	{x})=\alpha\in\mathbb{R}$.\smallskip\newline Let $y\in\dom s_{f}$ and let
	us assume, towards a contradiction, that $f_{l}(y)=+\infty$. Then for any
	sequence $\{y_{n}\}_{n}\subset\mathbb{R}^{d}$ that converges to $y$ we have
	$\underset{n\rightarrow\infty}{\lim\inf}\,f_n(y_{n})=+\infty.$ Note that the
	inequality
	\begin{equation}
		f_{n}(x)\,\geq\,f_{n}(y)-s_{f_{n}}(y)\,\Vert x-y\Vert\label{eq:ortega0}
	\end{equation}
	is valid for all $n\geq1$ and all $x,y\in\dom\,f_{n}$. Take the sequence
	$x_{n}\rightarrow\bar{x}$ given by \eqref{eq:Ortega_5000}, which verifies that ${\underset{n\rightarrow\infty}{\lim\inf
		}f_{n}(x_{n})=f_{l}(\bar{x})}$, and choose $\{y_{n}\}_{n}\subset\mathbb{R}^{d}$
	such that $\left(  y_{n},s_{f_{n}}(y_{n})\right)  \rightarrow(y,s_{f}(y))$
	(\textit{c.f.} Remark~\ref{rem: infsAttained}). Then replacing~$x$ by $x_{n}$
	and $y$ by $y_{n}$ in (\ref{eq:ortega0}) above, we easily deduce
	$f_{l}(\bar x)=+\infty$, which is a contradiction. \\
	Therefore, $f_{l}(y)<+\infty$.\smallskip\newline 
	On the other hand, setting $\sigma=\underset{n\rightarrow
		\infty}{\sup}s_{f_{n}}(x_{n})<+\infty,$ we easily see that for any sequence
	$\{y_{n}\}_{n}$ converging to $y$ we have:
	\[
	\liminf_{n\rightarrow\infty}\,f_{n}(y_{n})\,\geq\,\underset{n\rightarrow
		\infty}{\lim\inf}\,\left\{  \,f_{n}(x_{n})-s_{f_{n}}(x_{n})\,\Vert x_{n}
	-y_{n}\Vert\,\right\}  \,\geq\,\,f(\bar x)-\sigma\,\Vert \bar x-y\Vert\,>-\infty,
	\]
	yielding $f_{l}(y)\in\mathbb{R}$. Since $y$ is an arbitrary vector in $\dom
	s_{f}$, we have that $\dom s_{f}\subset\dom\,f_{l}$.\smallskip\newline Let us
	now show that $f_{u}(y)<+\infty$. Indeed, assuming the contrary, for any
	sequence $\{y_{n}\}_{n}\subset\mathbb{R}^{d}$ that converges to $y$, we would
	have $\underset{n\rightarrow\infty}{\limsup}f_{n}(y_{n})=+\infty$. Evoking
	again Remark~\ref{rem: infsAttained}, we can take $\{y_{n}\}_{n}$ such that
	$\left(  y_{n},s_{f_{n}}(y_{n})\right)  \rightarrow(y,s_{f}(y)).$ Therefore,
	for every $n\geq1$ we would have:
	\[
	f_{n}(x_{n})\,\geq\,f_{n}(y_{n})-s_{f_{n}}(y_{n})\,\Vert x_{n}-y_{n}\Vert
	\]
	and taking limsup at both sides of the above inequality, we would obtain
	$\alpha=f_{u}(\bar{x})=+\infty$, which is a contraction. Therefore,recalling that $-\infty < f_{l}(y)$ and that $f_{l}\leq f_u$, we
	deduce the inclusion
	\[
	\mathrm{dom}\,s_{f}\subset\mathrm{dom\,}f_{l}\cap\mathrm{dom\,}f_{u}.
	\]
	\newline Let us now show that $\mathrm{dom\,}f_{l}\cup\mathrm{dom\,}
	f_{u}\subset\overline{\mathrm{dom}}\,s_{f}.$ To this end, let $y\notin
	\overline{\mathrm{dom}}\,s_{f}$ and let $\varepsilon>0$ be such that
	$\overline{B}(y,\varepsilon)\subset\mathbb{R}^{d}\setminus\overline{\mathrm{dom}}\,s_{f}$.\smallskip\newline 
	We claim that $s_{f_{n}}\rightarrow\infty$ uniformly on $B(y,\varepsilon)$. Indeed, otherwise, there would exist $M>0$, a strictly
	increasing sequence $\{k_{n}\}_{n}\subset\mathbb{N}$ and a sequence
	$\{z_{n}\}_{n}\subset B(y,\varepsilon)$ such that ${s_{f_{k_{n}}}(z_{k_{n}})<M}$
	for all $n\in\mathbb{N}$. It follows that $s_{f}(z)\leq M$ for any cluster point $z\in
	\overline{B}(y,\varepsilon)$ of $\{z_{k_{n}}\}_{n}$ leading to a contradiction.\smallskip
	\newline We now set
	\[
	M_{n}:=\,\underset{z\in B(y,\varepsilon)}{\inf}\,s_{f_{n}}(z),\quad n\geq1,
	\]
	and observe that $M_{n}\rightarrow+\infty$ as $n\rightarrow\infty$. With this
	in mind, let us show that $f_{l}(y)=+\infty$.\smallskip\newline We proceed by
	contradiction: assume that $f_{l}(y)\in\mathbb{R}$, that is, for some sequence
	$\{y_{n}\}_{n}\subset B(y,\varepsilon)$ converging to $y$ we have
	$\underset{n\rightarrow\infty}{\liminf}\,f_{n}(y_{n})<\infty$. Then for $n$ sufficiently large (say $n\geq N$), let $\gamma
	_{n}:[0,\infty)\rightarrow\mathbb{R}^{d}$ be the steepest descend curve of the
	convex function $f_{n}$ starting at $y_{n},$ that is,
	\[
	\dot{\gamma}_{n}\in-\partial f_{n}(\gamma_{n})\qquad\text{and}\qquad\gamma
	_{n}(0)=y_{n}.
	\]
	Let further $\{t_{n}\}_{n}\subset(0,\infty)$ be the \textit{least escape-time}
	sequence defined by
	\[
	t_{n}=\,\inf\,\{t>0:~\gamma_{n}(t)\in\mathbb{R}^{d}\setminus B(y,\varepsilon
	)\}.
	\]
	In other words, $t_{n}>0$ is the first instant where the steepest descent
	curve $\gamma_{n}$ escapes from the ball $B(y,\varepsilon)$. Thus, $\gamma
	_{n}(t_{n})\in\partial B(y,\varepsilon)$ for all $n\geq N$. Since
	$\Vert\dot{\gamma}_{n}(\tau)\Vert=s_{f_{n}}(\gamma_{n}(\tau))$ (\textit{c.f.}
	\cite[Theorem~17.2.2]{ABM2014-book}) and since the length of the curve
	$\gamma_{n}$ in $[0,t_{n}]$ is larger than the distance $\mathrm{dist\,}(y_{n},\partial B(y,\varepsilon))=\varepsilon-\Vert y_{n}-y\Vert$ of the
	initial point $\gamma_{n}(0)=y_{n}$ to the boundary, we can write
	\[
	f_{n}(\gamma_{n}(t_{n}))\,=\,\,f_{n}(y_{n})-\int_{0}^{t_{n}}s_{f_{n}}
	(\gamma_{n}(\tau))\,\Vert\dot{\gamma}_{n}(\tau)\Vert\,d\tau\,\leq\,f_{n}
	(y_{n})-(\varepsilon-\Vert y_{n}-y\Vert)\,M_{n},
	\]
	concluding that $\underset{n\rightarrow\infty}{\liminf}f_{n}(\gamma_{n}
	(t_{n}))=-\infty$. However, convexity of $f_{n}$ at $x_{n}\in\dom\partial
	f_{n}$ yields that for all $n\in\mathbb{N}$ we have:
	\[
	f_{n}(\cdot)\geq g_{n}(\cdot):=f_{n}(x_{n})-s_{f_{n}}(x_{n})\,\Vert\cdot
	-x_{n}\Vert,
	\]
	and consequently,
	\begin{align*}
		\underset{n\rightarrow\infty}{\liminf}\,f_{n}(\gamma_{n}(t_{n})) &
		\geq\,\underset{n\rightarrow\infty}{\liminf}\,g_{n}(\gamma_{n}(t_{n}))\\
		&  \geq\,\underset{n\rightarrow\infty}{\lim}\,f_{n}(x_{n}
		)\,-\,\underset{n\rightarrow\infty}{\limsup}\,\left\{  s_{f_{n}}
		(x_{n})\,\left(  \Vert x_{n}-y\Vert+\Vert y-\gamma_{n}(t_{n})\Vert\right)
		\right\}  \\
		&  \geq\,f(x)-\left(  \underset{n\in\mathbb{N}}{\sup}\,s_{f_{n}}
		(x_{n})\right)  \left(  \Vert x-y\Vert+\varepsilon\right)  >-\infty,
	\end{align*}
	which is a contradiction.\smallskip\newline Therefore, $\mathrm{dom\,}
	f_{l}\subset\overline{\mathrm{dom}}\,s_{f}.$ Since $f_{l}\leq f_{u},$ the
	proof is complete. \hfill$\Box$
	
	
	\section{From slope convergence to epigraphical convergence}
	
	\label{sec: main}
	
	In this section we establish the difficult part of our main result, which
	states that up to a normalization condition, slope epigraphical convergence
	yields epigraphical convergence of the functions. This will be done in two
	stages: in Subsection~\ref{subsec: fu leq f} we show that $f_{u}\leq f$ while
	in Subsection~\ref{subsec: gap} we will control the gap between $f_{u}$ and
	$f_{l}$, then use~\eqref{eq:characEpiConvergence} to deduce our result.
	

	\subsection{Domination of the upper epigraphical limit}
	
	\label{subsec: fu leq f}
	
	We start with the following proposition:
	
	\begin{proposition}
		\label{Prop: f_u convex} Let $f_{n}:\mathbb{R}^{d}\rightarrow\mathbb{R}
		\cup\{+\infty\},$ $n\in\mathbb{N}$, be convex functions. Then
		\[
		f_{u}={\eLimsup}\,f_{n}
		\]
		is a convex lower semicontinuous function.
	\end{proposition}
	
	\noindent\textit{Proof.} We have already seen  in Remark~\ref{rem: infsAttained} that $f_u$ is lower semicontinuous.
	Let us now prove its convexity: to this end, consider any $x,y\in\mathbb{R}^{d}$ and $\lambda\in\lbrack0,1]$. 
	By Remark~\ref{rem: infsAttained} there exist sequences
	$\{x_{n}\}_{n},\{y_{n}\}_{n}\subset\mathbb{R}^{d}$ converging to $x$ and $y$
	respectively, such that
	\[
	\limsup_{n\rightarrow\infty}\,f_{n}(x_{n})=f_{u}(x)\qquad\text{and}
	\qquad\limsup_{n\rightarrow\infty}\,f_{n}(y_{n})=f_{u}(y).
	\]
	Using the fact that the functions $f_{n}$ are convex we have:
	\begin{align*}
		f_{u}(\lambda x+(1-\lambda)y)  &  =\inf_{z_{n}\rightarrow\lambda
			x+(1-\lambda)y}\limsup_{n\rightarrow\infty}\,f_{n}(z_{n})\,\leq\,\limsup
		_{n\rightarrow\infty}\,f_{n}(\lambda x_{n}+(1-\lambda)y_{n})\\
		&  \leq\,\limsup_{n\rightarrow\infty}\,\left\{  \lambda f_{n}(x_{n}
		)+(1-\lambda)f_{n}(y_{n})\,\right\}  \leq\,\lambda f_{u}(x)+(1-\lambda
		)f_{u}(y).
	\end{align*}
	This proves the convexity of $f_{u}$. \hfill$\Box$ 
	
	\bigskip
	
	Let us also recall (Remark~\ref{rem: infsAttained}) that the function
	$f_{l}={\eLiminf}\,f_{n}\ $is also lower semicontinuous.
	
	\begin{lemma}
		\label{lemma: sfu leq sf v1} Let $f,\{f_{n}\}_{n}:\mathbb{R}^{d}
		\rightarrow\mathbb{R}\cup\{+\infty\}$ be proper convex lower semicontinuous
		functions. Assume that $\{s_{f_{n}}\}_{n}$ epigraphically converges to $s_{f}$
		and that $f_{u}(\bar{x})\in\mathbb{R}$ for some $\bar{x}\in\mathrm{dom\,}
		s_{f}$. Then
		\[
		s_{f_{u}}\leq s_{f}.
		\]
		
	\end{lemma}
	
	\noindent\textit{Proof.} In view of Lemma~\ref{lemma: domInclusions} we have
	$\dom s_{f}\subset\dom\,f_{u}$. Let $y\in\dom s_{f}$ and let $\{y_{n}\}_{n}$
	be such that
	\[
	(y_{n},s_{f_{n}}(y_{n}))\rightarrow(y,s_{f}(y)).
	\]
	The sequence $\{y_{n}^{\ast}\}_{n}:=\{\partial^{\circ}f_{n}(y_{n})\}_{n}$ is
	then bounded. By Lemma~\ref{lemma: everySeqTight} we have that
	\[
	f_{u}(y)=\,\underset{n\rightarrow\infty}{\limsup}\,f_{n}(y_{n}).
	\]
	Passing to a subsequence, we may assume that for some $y^{\ast}\in
	\mathbb{R}^{d}$
	\[
	\underset{n\rightarrow\infty}{\lim}(y_{k_{n}},y_{k_{n}}^{\ast
	},f_{k_{n}}(y_{k_{n}}))=(y,y^{\ast},f_u(y)).
	\]
	Since $s_{f_{n}}(y_{n})\rightarrow s_{f}(y),$ it follows easily that $\Vert
	y^{\ast}\Vert=s_{f}(y)$. Furthermore, for any $z\in\mathbb{R}^{d}$ and any
	sequence $\{z_{n}\}_{n}$ converging to $z$ we have
	\begin{align*}
		\limsup_{n\rightarrow\infty}\,f_{n}(z_{n})\, &  \geq\,\limsup_{n\rightarrow
			\infty}\,f_{k_{n}}(z_{k_{n}})\\
		&  \geq\,\limsup_{n\rightarrow\infty}\,\left\{  f_{k_{n}}(y_{k_{n}})+\langle
		y_{k_{n}}^{\ast},z_{k_{n}}-y_{k_{n}}\rangle\right\}  \\
		&  =\,f_{u}(y)+\langle y^{\ast},z-y\rangle.
	\end{align*}
	Since $\{z_{n}\}_{n}$ is an arbitrary sequence, we deduce $f_{u}(z)\geq
	f_{u}(y)+\langle y^{\ast},z-y\rangle.$ Since $z$ is arbitrary, we obtain that
	$y^{\ast}\in\partial f_{u}(y)$. Thus, $s_{f_{u}}(y)\leq\Vert y^{\ast}
	\Vert=s_{f}(y)$.\smallskip\newline If $y\notin\mathrm{dom}\,s_{f},$ the
	inequality $s_{f_{u}}(y)\leq s_{f}(y)\equiv+\infty$ is obvious. The proof is
	complete. \hfill$\Box$
	
	\bigskip
	
	The above lemma will be used in the following result.
	
	\begin{proposition}
		\label{prop: fu leq f v2} Let $f,\{f_{n}\}_{n}:\mathbb{R}^{d}\rightarrow
		\mathbb{R}\cup\{+\infty\}$ be proper convex lower semicontinuous functions.
		Assume that $\{s_{f_{n}}\}_{n}$ converges epigraphically to $s_{f}$ and that
		$\inf\,f=\inf f_{u}\in\mathbb{R}$. Then
		\[
		f_{u}\leq f.
		\]
		
	\end{proposition}
	
	\noindent\textit{Proof}.
	Let $\bar{x}\in\mathrm{dom}\,s_{f}$ and let us show that $f_{u}(\bar{x}
	)\in\mathbb{R}$. Indeed, one obviously has $f_{u}(\bar{x})\geq\inf
	\,f_{u}>-\infty$. Reasoning towards a contradiction, let us assume that
	$f_{u}(\bar{x})=+\infty$. Then for any sequence $\{x_{n}\}_{n}$ converging to
	$\bar{x}$, we have that $\underset{n\rightarrow\infty}{\limsup}\,f_{n}
	(x_{n})=+\infty$. We may choose $\{x_{n}\}_{n}$ so that $(x_{n},s_{f_{n}
	}(x_{n}))\rightarrow(\bar{x},s_{f}(\bar{x}))$, ensuring in particular that
	$\{s_{f_{n}}(x_{n})\}_{n}$ is bounded. Then for every $y\in\mathbb{R}^{d}$ and
	every sequence $\{y_{n}\}_{n}\subset\mathbb{R}^{d}$ such that $y_{n}
	\rightarrow y$ we have:
	\begin{align*}
		\limsup_{n\rightarrow\infty}\,f_{n}(y)\, &  \geq\,\limsup_{n\rightarrow\infty
		}\,\left\{  f_{n}(x_{n})-s_{f_{n}}(x_{n})\,\Vert x_{n}-y_{n}\Vert\right\}  \\
		&  =\,\limsup_{n\rightarrow\infty}\,f_{n}(x_{n})-s_{f}(\bar{x})\,\Vert\bar
		{x}-y\Vert\,=\,+\infty.
	\end{align*}
	This yields that $f_{u}\equiv+\infty$, which is a contradiction.\smallskip
	\newline Therefore, $f_{u}(\bar{x})\in\mathbb{R}$ and we can apply
	Lemma~\ref{lemma: sfu leq sf v1} to get that $s_{f_{u}}\leq s_{f}$. The
	conclusion follows from Theorem~\ref{theo: comparisonPrinciple} (comparison
	principle). \hfill$\Box$
	
	\bigskip
	
	Forthcoming Lemma~\ref{lem: curve convergence} provides a criterium for a
	limit of steepest descent curves (of convex functions converging
	epigraphically to a limit function) to be a steepest descent curve of the limit
	function. This is an intermediate result, which will be further refined in
	Lemma~\ref{lem: f_u against f} and eventually lead to
	Proposition~\ref{prop: f_u leq f v1} (domination of $f_{u}$ by $f$).
	\smallskip\newline We shall first need the following result.
	
	\begin{proposition}
		\label{Prop: curves} Let $\{\gamma_{n}\}_{n}$ be a sequence of Lipschitz
		curves from $[0,+\infty)$ to $\mathbb{R}^{d}$. Assume that the sequence
		$\{\gamma_{n}(0)\}_{n}$ is bounded and that all Lipschitz constants of the
		curves $\{\gamma_{n}\}_{n}$ are bounded by a constant $K>0$, that is,
		$\mathrm{Lip}(\gamma_{n})\leq K$ for all $n\geq1.$ Then, there exists an
		increasing sequence $\{k(n)\}_{n}$ such that $\{\gamma_{k(n)}\}_{n}$ converges
		uniformly on compact sets to a Lipschitz curve $\gamma:[0,+\infty)\rightarrow\mathbb{R}^{d}$
		and the sequence of its tangents $\{\dot{\gamma}_{k(n)}|_{[0,T]}\}_{n}$
		converges weakly to $\dot{\gamma}|_{[0,T]}$ in $\mathcal{L}^{2}
		([0,T];\mathbb{R}^{d})$, for any $T>0$.
	\end{proposition}
	
	\noindent\textit{Proof.} Let us first assume that the curves $\{\gamma
	_{n}\}_{n}$ are defined on $[0,1]$. Since $\{\gamma_{n}(0)\}_{n}$ is
	relatively compact and $\mathrm{Lip}(\gamma_{n})\leq K$ for all $n\in
	\mathbb{N}$, we can apply Arzel\`{a}--Ascoli theorem to get a subsequence
	$\{\gamma_{k(n)}\}_{n}$ which converges uniformly to some continuous curve
	$\gamma$ on $[0,1]$. It follows easily that $\gamma$ is Lipschitz with
	$\mathrm{Lip}(\gamma_{n})\leq K$. Since $\{\dot{\gamma}_{k(n)}\}_{n}$ is
	bounded on $\mathcal{L}^{2}([0,T];\mathbb{R}^{d})$ (in fact $\Vert\dot{\gamma
	}_{k(n)}\Vert_{\mathcal{L}^{2}}\leq K$), by the Eberlein--\u{S}mulian theorem,
	there exists a subsequence $\{k^{\prime}(k(n))\}_{n}$ which we denote by
	$\{\bar{k}(n)\}_{n}$ (that is, $\bar{k}=k^{\prime}\circ k$), such that
	$\{\dot{\gamma}_{\bar{k}(n)}\}_{n}$ converges weakly to $\nu:[0,1]\rightarrow
	\mathbb{R}^{d}$. Notice that, for each $n\geq1$ we have that
	\[
	\gamma_{\bar{k}(n)}(t)=\gamma_{\bar{k}(n)}(0)+\int_{0}^{t}\dot{\gamma}
	_{\bar{k}(n)}(s)ds,~\text{for all }t\in\lbrack0,1].
	\]
	Taking the limit as $n\rightarrow+\infty$, we obtain
	\[
	\gamma(t)=\gamma(0)+\int_{0}^{t}\nu(s)ds,~\text{for all }t\in\lbrack0,1].
	\]
	Therefore, $\dot{\gamma}(t)=\nu(t),$ for all $t\in\lbrack0,1]$ and the
	assertion follows.\smallskip\newline The general case follows easily: if the
	curves $\{\gamma_{n}\}_{n}$ are defined on $[0,+\infty)$, we fix $T>0$ and
	proceed as before for the restricted curves $\{\gamma_{n}|_{[0,T]}\}_{n}$. The
	result follows via a standard diagonal argument. \hfill$\Box$ \bigskip
	
	We are now ready to prove our lemma.
	
	\begin{lemma}
		\label{lem: curve convergence}Let $f,\{f_{n}\}_{n}:\mathbb{R}^{d} \rightarrow\mathbb{R}\cup\{+\infty\}$ be proper convex lower semicontinuous
		functions such that
		\[
		s_{f_{n}}\xrightarrow{e}s_{f}\qquad\text{and}\qquad\inf f>-\infty.
		\]
		Assume that there is a sequence $(x_{n},x_{n}^{\ast},f_{n}(x_{n}))\in\triangle f_{n},$ $n\geq1$ such that
		\[
		\lim_{n\rightarrow\infty}\,(x_{n},x_{n}^{\ast},f_{n}(x_{n}))=(\bar{x},\bar
		{x}^{\ast},f(\bar{x}))\in\triangle f.
		\]
		For every $n\geq1,$ let $\gamma_{n}:[0,+\infty)\rightarrow\mathbb{R}^{d}$ denote the steepest descent curve of $f_{n}$
		starting from the point $x_{n}\in\dom\,f_{n},$ 
		and let us assume that: \newline
		\phantom{jo}$\mathrm{(i).}$ $\{\gamma_{n}\}_{n}$ converges to some Lipschitz curve~$\nu$ uniformly on compact sets; and \\
		\phantom{jo}$\mathrm{(ii).}$ for all $T>0,$ the tangents $\{\dot{\gamma}_{n}|_{[0,T]}\}_{n}$ converge to $\dot{\nu}|_{[0,T]}$
		weakly on $\mathcal{L}^{2}([0,T],\mathbb{R}^{d})$. \\
		Then $\nu$ is a steepest descent curve for the function $f_{u}$.
	\end{lemma}
	
	\noindent\textit{Proof.} Thanks to Lemma~\ref{lemma: domInclusions}, we know
	that $\mathrm{dom}\,s_{f}\subset\mathrm{dom\,}f_{u}\,\subset\,\overline
	{\mathrm{dom}}\,s_{f},$ which yields
	\[
	\ri(\dom s_{f})=\ri(\dom\,f_{u}).
	\]
	Since $\gamma_{n}$ is a steepest descent curve of $f_{n}$ emanating from
	$x_{n}$, for every $t>0$ and $n\in\mathbb{N}$ we have:
	\[
	f_{n}(\gamma_{n}(t))\leq f_{n}(x_{n})\quad\text{ and }\quad s_{f_{n}}
	(\gamma_{n}(t))\leq s_{f_{n}}(x_{n}).
	\]
	It follows easily from our hypothesis that the sequence $\{s_{f_{n}}
	(\gamma_{n}(t))\}_{n}$ is bounded. Since ${\gamma_{n}(t)\rightarrow\nu(t)}$,
	Lemma~\ref{lemma: everySeqTight} entails that
	\[
	f_{u}(\nu(t))=\,\limsup_{n\rightarrow\infty}\,f_{n}(\gamma_{n}(t))\,\leq
	\,\limsup_{n\rightarrow\infty}\,f_{n}(x_{n})\,=\,f(\bar{x})\,<\,+\infty.
	\]
	Thus, for every $t>0$, we have $\nu(t)\in\dom\,f_{u}$.\smallskip\newline Let
	$y\in\ri(\dom\,f_{u})$ and let $(y_{n},y_{n}^{\ast})\in\partial f_{n}$ be such
	that $\{y_{n}\}_{n}$ converges to $y$ and the sequence $\{\Vert y_{n}^{\ast
	}\Vert\}_{n}=\{s_{f_{n}}(y_{n})\}_{n}$ converges to $s_{f}(y)$. Passing to a
	subsequence $\{(y_{k_{n}},y_{k_{n}}^{\ast})\}_{n}$, we obtain
	\[
	f_{u}(y)=\,\underset{n\rightarrow\infty}{\lim}f_{k_{n}}(y_{k_{n}}
	)\quad\text{and}\quad\underset{n\rightarrow\infty}{\lim}y_{k_{n}}^{\ast
	}=y^{\ast}\quad\text{for some }y^{\ast}\in\mathbb{R}^{d}.
	\]
	Using the same argument as in the proof of Lemma~\ref{lemma: sfu leq sf v1} we
	deduce that $y^{\ast}\in\partial f_{u}(y)$. Then, for any bounded Borel set
	$A\subset\lbrack0, +\infty)$ and any $n\in\mathbb{N}$ we have
	\[
	0\leq\int_{A}\left\langle \,y_{k_{n}}^{\ast}+\dot{\gamma}_{k_{n}
	}(t),\,y_{k_{n}}-\gamma_{k_{n}}(t)\,\right\rangle \,dt.
	\]
	Taking the limit as $n\rightarrow\infty$ we obtain (thanks to our assumption) that
	\[
	0\leq\int_{A}\left\langle y^{\ast}+\dot{\nu}(t),\,y-\nu(t)\right\rangle \,dt.
	\]
	Since $A\subset\lbrack0,+\infty)$ is an arbitrary bounded Borel set, we deduce
	that $\langle y^{\ast}+\dot{\nu}(t),y-\nu(t)\rangle\geq0$ for a.e.
	$t\in\lbrack0,+\infty)$. Since $y\in\ri(\dom\,f_{u})$ is arbitrary, we can
	take a sequence $\{(z_{n},z_{n}^{\ast})\}_{n}\subset\partial f_{u}$ such
	that $\Vert z_{n}^{\ast}\Vert=s_{f}(z_{n}),$ for all $n\in\mathbb{N}$, and
	$\{z_{n}\}_{n}$ is dense in $\ri(\dom\,f_{u})$, obtaining
	\[
	0\leq\langle z_{n}^{\ast}+\dot{\nu}(t),z_{n}-\nu(t)\rangle,~\forall_{a.e.}
	t\in\lbrack0,+\infty),~\forall n\in\mathbb{N}.
	\]
	Thus, applying Proposition~\ref{Prop: subdiff r-int}, we deduce that
	\[
	\dot{\nu}(t)\in-\partial f_{u}(\nu(t)),~\forall_{a.e.}t\in\lbrack0,+\infty).
	\]
	The proof is complete. \hfill$\Box$
	
	\bigskip
	Before we proceed, let us recall an important technical result ensuring that absolutely continuous curves verifying an integrability condition for the slope of a convex function must be infimizing. This result is essentially known, but we include a proof for completeness, since the precise statement that we use below in not directly available in the literature. A~strengthened version (which is also contained in the proposition below) can be obtained if the curve is a steepest descent curve of another function, see \cite[Lemma~3.1]{PSV2021}. A discretized version has been used in \cite{TZ2023}.
	
	\begin{proposition}
		[infimizing curves by integrability of slope]\label{lem_salas} Let
		$g:\mathbb{R}^{d}\rightarrow\mathbb{R}\cup\{+\infty\}$ be a proper convex
		lower semicontinuous function, and let $\gamma:[0,+\infty)\mapsto \mathbb{R}^d$ be an absolutely continuous curve. The following assertions hold:
		\begin{enumerate}
			\item[$\mathrm{(i).}$] If $\gamma$ satisfies that
			\begin{equation}\label{eq:IntegrabilitySlope}
				\liminf s_g(\gamma(t)) = 0\qquad\text{ and} \qquad \int_0^{+\infty} s_g(\gamma(t))\,\|\dot{\gamma}(t)\|\,dt < +\infty,
			\end{equation}
			then $$\underset{t\to +\infty}{\liminf} g(\gamma(t)) = \inf g.$$   
			\item[$\mathrm{(ii).}$] Let $f$ be another proper convex lower semicontinuous function such that 
			\[
			s_{f}(x)\geq s_{g}(x),\quad\text{for all }x\in\mathbb{R}^{d}\qquad
			\text{and\qquad}\inf\,f>-\infty.
			\]
			If $\gamma:[0,+\infty)\mapsto \mathbb{R}^d$ is the steepest descent curve for $f$ starting at a point
			$\bar{x}\in\mathrm{dom\,}f$, then
			\[
			\underset{t\rightarrow +\infty}{\lim}\,g(\gamma_{\bar{x}}(t))=\inf\,g.
			\]
			Moreover, if $f(\bar{x}) = g(\bar{x})$, then $\inf g \geq \inf f$.
		\end{enumerate}
	\end{proposition}
	
	\noindent \textit{Proof}. (i). Integrability of the slope implies that $\gamma (t)\in \dom\partial g$ for almost every $t\geq 0$. 
	Thus, for almost every $t\geq 0$, we can define $h(t)=\partial ^{\circ }f(\gamma (t))$ so
	that $\Vert h(t)\Vert =s_{f}(\gamma (t))$. Then via the standard
	subdifferential calculus (chain rule) for convex functions (see, e.g., \cite[Proposition~17.2.5]{ABM2014-book}) we deduce: 
	\begin{equation*}
		g(\gamma (t))-g(\gamma (0))=\int_{0}^{t}\langle h(\tau ),\dot{\gamma}(\tau
		)\rangle d\tau \leq \int_{0}^{+\infty }s_{g}(\gamma (t))\Vert \dot{\gamma}
		(t)\Vert dt<+\infty .
	\end{equation*}
	Thus, $\underset{t\rightarrow +\infty }{\liminf }\, g(\gamma (t))\leq \,
	\underset{t\rightarrow +\infty }{\limsup }\,g(\gamma (t))<+\infty $.\smallskip 
	
	\textit{Claim.} There exists an increasing sequence $\{t_{n}\}_{n\in \mathbb{N}}$ such that
	\begin{equation}
		t_{n}\nearrow +\infty \qquad \text{and}\qquad s_{g}(\gamma (t_{n}))=\,\underset{t\in \lbrack t_{0},t_{n}]}{\min }\,s_{g}(\gamma (t)).
		\label{eq:salov}
	\end{equation}
	\smallskip \newline
	By assumption we have that $\underset{t\rightarrow +\infty }{\liminf }
	\,s_{g}(\gamma (t))=0.$ If $s_{g}(\gamma (t))=0$ recurrently as 
	$t\rightarrow +\infty $, then we choose $\{t_{n}\}_{n}$ to be any increasing
	sequence with $t_{n}\rightarrow +\infty $ and $s_{g}(\gamma (t_{n}))=0$ for
	every $n\in \mathbb{N}$. Otherwise, we fix $t_{0}\geq \sup \{t\in \lbrack
	0,+\infty )\ :\ s_{g}(\gamma (t))=0\}+1$ and we define:
	\begin{equation*}
		M_{n}:=\underset{t\in \lbrack t_{0},t_{0}+n]}{\argmin}s_{g}(\gamma (\cdot
		))\qquad \text{and}\qquad t_{n}=\max M_{n}.
	\end{equation*}
	Since $s_{g}$ is lower semicontinuous by convexity of $g$ (see, e.g., \cite{AGS-2008}), the sequence $\{t_{n}\}_{n}$ is well defined and \eqref{eq:salov}  holds.
	
	\bigskip 
	
	Now, take any $v\in \dom g$ and $T_{\varepsilon }\geq t_{0}$ large enough
	such that 
	\begin{equation*}
		\int_{T_{\varepsilon }}^{+\infty }s_{g}(\gamma (t))\Vert \dot{\gamma}
		(t)\Vert dt\leq \varepsilon .
	\end{equation*}
	Using convexity, Cauchy-Schwarz inequality and (\ref{eq:salov}) we deduce
	that for all $t_{n}>T_{\varepsilon }$ we have: 
	\begin{align*}
		g(\gamma (t_{n}))& \leq g(v)+\,\langle h(t_{n}),\gamma (t_{n})-v\rangle  \\
		& \leq g(v)+\,|\langle h(t_{n}),v\rangle |\,+\,|\langle h(t_{n}),\gamma
		(t_{\varepsilon })\rangle |\,+\,\int_{T_{\varepsilon }}^{t_{n}}|\langle
		h(t_{n}),\dot{\gamma}(s)\rangle |ds \\
		& \leq g(v)+\,|\langle h(t_{n}),v\rangle |\,+\,|\langle h(t_{n}),\gamma
		(t_{\varepsilon })\rangle |\,+\,\int_{T_{\varepsilon }}^{t_{n}}s_{g}(\gamma
		(t_{n}))\,\Vert \dot{\gamma}(s)\Vert ds \\
		& \leq g(v)+\,s_{g}(\gamma (t_{n}))\Vert v\Vert +\,s_{g}(\gamma
		(t_{n}))\Vert \gamma (t_{\varepsilon })\Vert +\,\int_{T_{\varepsilon
		}}^{t_{n}}s_{g}(\gamma (s))\Vert \dot{\gamma}(s)\Vert ds \\
		& \xrightarrow{n\to\infty}g(v)+\int_{T_{\varepsilon }}^{+\infty
		}s_{g}(\gamma (s))\Vert \dot{\gamma}(s)\Vert ds\leq g(v)+\varepsilon .
	\end{align*}
	Thus, for every $\varepsilon >0$ and every $v\in \dom g$, we have that 
	\begin{equation*}
		\liminf_{t\rightarrow +\infty }g(\gamma (t))\leq \liminf_{n\rightarrow
			\infty }g(\gamma (t_{n}))\leq g(v)+\varepsilon .
	\end{equation*}
	Thus, $\underset{t\rightarrow +\infty }{\liminf }g(\gamma (t))=\inf g$,
	finishing the proof of this part.\medskip \newline
	(ii). The first conclusion of the second part is given by \cite[Lemma 3.1]{PSV2021} and the proof is very similar of the latter development, but using that $s_f(\gamma(t))$ is nonincreasing as $t\to +\infty$. For the last part, it is enough to write
	
	\begin{align*}
		\inf g-g(\bar{x}) 
		&  =\liminf_{t\rightarrow\infty}\int_{0}^{t}\frac{d}{dt}[g\circ\gamma](\tau)\,d\tau = \liminf_{t\rightarrow\infty}\int_{0}^{t} \langle\partial^{\circ}g(\gamma(\tau)),\dot{\gamma}(\tau)\rangle\,d\tau \\
		& \underbrace{ \geq}_{Cauchy-Schwarz}-\limsup_{t\rightarrow\infty}\int_{0}^{t}s_{g}(\gamma(\tau))\|\dot{\gamma}(\tau)\|\,d\tau\,
		\underbrace{\geq}_{s_f\geq s_g}-\lim_{t\rightarrow\infty}\int_{0}^{t}s_{f}(\gamma(\tau ))^{2}\,d\tau\, \\
		& =\lim_{t\rightarrow\infty}\int_{0}^{t}\langle\partial^{\circ}f(\gamma(\tau)),\dot{\gamma}(\tau)\rangle\,d\tau =\lim_{t\rightarrow\infty}\int_{0}^{t}\frac{d}{dt}[f\circ\gamma](\tau))\,d\tau\,=\inf f\,-f(\bar{x}).
	\end{align*}
	The result follows. \hfill$\Box$

	\bigskip
	
	We are now ready to obtain an enhanced version of
	Lemma~\ref{lem: curve convergence}.
	
	\begin{lemma}
		\label{lem: f_u against f} Under the same assumptions as in
		Lemma~\ref{lem: curve convergence} we conclude:
		\[
		\inf f_{u}=\inf f\qquad\text{and}\qquad f_{u}\leq f.
		\]
		
	\end{lemma}
	
	\noindent\textit{Proof.} By hypothesis $\inf\,f>-\infty.$ Moreover, by
	Lemma~\ref{lemma: everySeqTight} we obtain $f_{u}(\bar{x})=f(\bar{x})$, while
	by Lemma~\ref{lem: curve convergence} (and following notation therein) the
	limit curve $\nu=\underset{n\rightarrow\infty}{\lim}\gamma_{n}$ is a steepest
	descent curve for (the convex function) $f_{u}$ starting at $\bar{x}=\nu(0)$,
	that is,
	\begin{equation}
		\forall_{\text{a.e.\thinspace}}t\in\lbrack0,+\infty):\;\dot{\nu}
		(t)=-\partial^{\circ}f_{u}(\nu(t))\,\text{,}\quad\Vert\dot{\nu}(t)\Vert
		=s_{f_{u}}(\nu(t))\quad\text{and}\quad f_{u}(\nu(t))\rightarrow\inf
		\,f_{u}.\label{eq:ortega-1}
	\end{equation}
	Let us also recall from Lemma~\ref{lemma: sfu leq sf v1} that
	\begin{equation}
		s_{f}(x)\geq s_{f_{u}}(x),\qquad\text{for all }x\in\mathbb{R}^{d}
		.\label{eq:ortega-2}
	\end{equation}
	
	By Proposition \ref{lem_salas}.(ii) it holds:
	\begin{equation}
		\inf f_{u}\geq\inf f>-\infty.\label{eq:claim}
	\end{equation}
	\smallskip\newline
	
	Let us set $M=\sup\{\Vert x_{n}^{\ast}\Vert:~n\in\mathbb{N}\} = \sup\{s_{f_n}(x_n)\, :~n\in\mathbb{N}\}$ and notice that for all $t\geq0$ and
	$n\geq1$ we have $s_{f_{n}}(\gamma_{n}(t))\leq~s_{f_{n}}(x_{n})\leq~M$. We deduce easily from Lemma~\ref{lemma: everySeqTight} that
	\[
	f_{u}(\nu(t))=\limsup_{n\rightarrow+\infty}\,f_{n}(\gamma_{n}(t)),\quad
	\text{for all\ }t\geq0.
	\]
	By Fatou's Lemma and (\ref{eq:ortega-2}) we have
	\begin{align}
		f_{u}(\nu(t)) &  =\limsup_{n\rightarrow\infty}\,f_{n}(\gamma_{n}
		(t))=\limsup_{n\rightarrow\infty}\,\left\{  f_{n}(\gamma_{n}(0))-\int_{0}
		^{t}s_{f_{n}}(\gamma_{n}(\tau))^{2}\,d\tau\right\}  \label{eq:jaime}\\
		&  =\,f(\bar{x})-\liminf_{n\rightarrow\infty}\,\int_{0}^{t}s_{f_{n}}
		(\gamma_{n}(\tau))^{2}\,d\tau\,\leq\,f(\bar{x})-\int_{0}^{t}\liminf
		_{n\rightarrow\infty}\,s_{f_{n}}(\gamma_{n}(\tau))^{2}\,d\tau\nonumber\\
		&  \leq\,f(\bar{x})-\int_{0}^{t}s_{f}(\nu(\tau))^{2}\,d\tau\,\leq
		\,f(\nu(t)).\nonumber
	\end{align}
	Therefore, we deduce:
	\[
	\int_{0}^{t}s_{f}(\nu(\tau))^{2}\,d\tau\,\leq f(\bar{x})-\inf f_{u}
	<+\infty,\qquad\text{for every\ }t\geq0
	\]
	and
	\[
	\inf\,f_{u}\,=\,\inf\,\left(  f_{u}\circ\nu\right)  \,\leq\,\liminf
	_{t\rightarrow\infty}\,\left(  f\circ\nu\right).
	\]
	Applying Proposition~\ref{lem_salas}.(i) to $g=f$ and $\gamma =\nu$, and noting that 
	\[
	\liminf s_{f_u}(\nu(t)) \leq \liminf s_{f}(\nu(t)) = 0,
	\]
	we get that $\liminf
	_{t\rightarrow\infty}\,\left(  f\circ\nu\right) = \inf f$. We conclude that
	\[
	\inf f_{u}=\inf f\in\mathbb{R}.
	\]
	The result follows by applying Proposition~\ref{prop: fu leq f v2}.\hfill
	$\Box$
	
	\bigskip
	
	We finish this subsection with the following proposition that, together with Proposition \ref{prop: fu leq f v2}, provides a partial result towards  our main theorem: If (ii) or (iii)  of Theorem \ref{theo: main} hold, then $f_u\leq f$.
	
	\begin{proposition}
		\label{prop: f_u leq f v1} Let $f,\{f_{n}\}_{n}:\mathbb{R}^{d}\rightarrow
		\mathbb{R}\cup\{+\infty\}$ be proper convex lower semicontinuous functions
		such that
		\[
		s_{f_{n}}\xrightarrow{e}s_{f}\qquad\text{and}\qquad\inf f>-\infty.
		\]
		Assume that there is a sequence
		\[
		(x_{n},x_{n}^{\ast},f_{n}(x_{n}))\in\triangle f_{n}\qquad\text{and}\qquad
		\lim_{n\rightarrow\infty}\,(x_{n},x_{n}^{\ast},f_{n}(x_{n}))=(\bar{x},\bar
		{x}^{\ast},f(\bar{x}))\in\triangle f.
		\]
		Then
		\[
		\inf\,f_{u}=\inf f\qquad\text{and}\qquad f_{u}\leq f.
		\]
		
	\end{proposition}
	
	\noindent\textit{Proof.} Let $\gamma_{n}$ be the steepest descent curve of
	$f_{n}$ starting at $x_{n}=\gamma_{n}(0)$. Set
	\[
	M=\,\underset{n\geq1}{\sup}\,\Vert x_{n}^{\ast}\Vert
	\]
	so that
	\[
	s_{f_{n}}(\gamma_{n}(t))\leq s_{f_{n}}(\gamma_{n}(0))=\Vert x_{n}^{\ast}
	\Vert\leq M\qquad\text{for all }t\geq0\text{ and }n\geq1.
	\]
	By a standard application of Arzel\`{a}-Ascoli theorem, for every strictly
	increasing sequence $\{k_{1}(n)\}_{n}$ there exists a subsequence
	$\{(k_{2}\circ k_{1})(n)\}_{n}$ that we simply denote by $\{k_{n}\}_{n}$ such
	that $\{\gamma_{k_{n}}\}_{n}$ uniformly converges to some Lipschitz curve on
	$[0,T]$, for every $T>0$ (as in the statement of Lemma~\ref{lem: curve convergence}).
	Up to a new subsequence, which we keep denoting as before, $\{\gamma_{k_{n}}\}_{n}$ converges to a Lipschitz curve~$\nu$ uniformly on bounded sets 
	and $\{\dot{\gamma}_{k_{n}} |_{[0,T]}\}_{n}$ converges weakly to $\dot{\nu}|_{[0,T]}$ in $\mathcal{L}^{2}([0,T],\mathbb{R}^{d})$.
	Let
	\[
	f_{u,k_{n}}:={{\eLimsup}}\,f_{k_{n}},\quad\text{for all }n\geq1.
	\]
	Thanks to Lemma~\ref{lem: f_u against f}, we have $f_{u,k_{n}}\leq f$. Since
	this holds true for any sequence $\{k_{n}\}_{n}$ such that $\{\gamma_{k_{n}
	}\}_{n}$ converges (as in the statement of Lemma~\ref{lem: f_u against f}), we
	can claim that $f_{u}\leq f$.\smallskip\newline Indeed, for any $y\in\dom
	s_{f}$, there is a sequence $\{y_{n}\}_{n}$ such that $(y_{n},s_{f_{n}}
	(y_{n}))\rightarrow(y,s_{f}(y))$. By Lemma~\ref{lemma: everySeqTight}, there
	exists a subsequence $\{k_{1}(n)\}_{n}$ such that $f_{u}(x)=\lim_{n}
	f_{k_{1}(n)}(y_{k_{1}(n)})$. By Arzel\`{a}-Ascoli theorem, there exists a
	sub-subsequence $\{k_{2}(k_{1}(n)\}\}_{n}$ such that for $k=k_{2}\circ k_{1}$
	the sequence of the steepest descend curves $\{\gamma_{k_{n}}\}_{n}$ converges
	to a curve $\nu$ (as in the statement of Lemma~\ref{lem: f_u against f}) and we get $f_{u,k_{n}}\leq f$. 
	Therefore, thanks to Lemma~\ref{lemma: everySeqTight}, we infer that
	\[
	f_{u}(y)=f_{u,k_{n}}(y)\leq f(y).
	\]
	Since $y$ is an arbitrary vector in $\dom s_{f}$, we obtain $f_{u}\leq f$ on
	$\dom s_{f}$. Now, recalling Proposition~\ref{Prop: f_u convex}, $f$ and
	$f_{u}$ are convex lower semicontinuous functions and it is enough to apply
	Proposition~\ref{prop: extension to the adherence} and
	Lemma~\ref{lemma: domInclusions} to conclude that $f_{u}\leq f$ on
	$\mathbb{R}^{d}$. \hfill$\Box$
	

	\subsection{Controlling the gap between upper and lower epigraphical limits}
	
	\label{subsec: gap}
	
	Let us first recall the following important result from \cite[Lemma~2.4]{CT1998}.
	
	\begin{proposition}\label{prop: Thibault Combari} Let $f_{n}:\mathbb{R}
		^{d}\rightarrow\mathbb{R}\cup\{+\infty\}$ be convex lower semicontinuous
		functions such that there exists a sequence $(x_{n},x_{n}^{\ast},f_{n}
		(x_{n}))\in\triangle f_{n},$ $n\geq1,$ such that
		\[
		\lim_{n\rightarrow\infty}\,(x_{n},x_{n}^{\ast},f_{n}(x_{n}))=(\bar{x},\bar
		{x}^{\ast},\alpha)\in\mathbb{R}^{d}\times\mathbb{R}^{d}\times\mathbb{R}.
		\]
		Then,
		\[
		\bar{x}^{\ast}\in\partial f_{u}(\bar{x})\cap\partial f_{l}(\bar{x}
		)\qquad\text{and}\qquad\alpha=f_{l}(\bar{x})=f_{u}(\bar{x}).
		\]
		
	\end{proposition}
	
	\noindent\textit{Proof.} By Lemma~\ref{lemma: everySeqTight}, we have
	$f_{l}(\bar{x})=f_{u}(\bar{x})=\lim_{n}f_{n}(\bar{x})$. Then, for any
	$y\in\mathbb{R}^{d}$ and any sequence $\{y_{n}\}_{n}\subset\mathbb{R}^{d}$
	converging to $y$, we have
	\[
	f_{n}(y_{n})\geq f_{n}(x_{n})+\langle x_{n}^{\ast},y_{n}-x_{n}\rangle
	,\text{\qquad for all }n\geq1.
	\]
	The desired conclusion follows by taking $\lim\sup$ and $\lim\inf$ to the
	above expression. \hfill$\Box$
	
	\bigskip
	
	The following result states that epigraphical convergence of the sequence of
	slope functions guarantees the local Lipschitz continuity of the lower
	epigraphical limit function $f_{l}$ under a mild condition.
	
	\begin{proposition}
		\label{prop: loc lipschitz f_l} Let $f,\{f_{n}\}_{n}:\mathbb{R}^{d}
		\rightarrow\mathbb{R}\cup\{+\infty\}$ be proper convex lower semicontinuous
		functions. Assume that $\{s_{f_{n}}\}_{n}$ epigraphically converges to $s_{f}
		$. Assume further that there is a sequence $\{x_{n}\}_{n}\subset\mathbb{R}
		^{d}$ converging to $\bar{x}$ such that $\{s_{f_{n}}(x_{n})\}_{n}$ is bounded
		and $\{f_{n}(x_{n})\}_{n}$ converges. Then, $f_{l}$ is locally Lipschitz on
		$\ri(\dom\,f_{l})$.
	\end{proposition}
	
	\noindent\textit{Proof.} Since $f$ is convex and lower semicontinuous, we know
	by Proposition~\ref{prop: convex rint} that $\ri(\dom s_{f})$ is a convex set.
	Thus, thanks to Lemma~\ref{lemma: domInclusions}, we have that
	$\ri(\dom\,f_{l})=\ri(\dom s_{f})$. Let $y,z\in\dom s_{f}$ and let
	$\{y_{n}\}_{n},\{z_{n}\}_{n}\subset\mathbb{R}^{d}$ be two sequences convergent
	to $y$ and $z$, such that $\{s_{f_{n}}(y_{n})\}_{n}$ and $(s_{f_{n}}
	(z_{n}))_{n}$ converge to $s_{f}(y)$ and $s_{f}(z)$, respectively. By
	Lemma~\ref{lemma: everySeqTight}, we get that $f_{l}(y)=\liminf_{n}f_{n}
	(y_{n})$ and $f_{l}(z)=\liminf_{n}\,f_{n}(z_{n})$. Take a subsequence
	$(k_{n})_{k}$ such that $f_{k_{n}}(y_{k_{n}})\rightarrow f_{l}(y)$. Then,
	\begin{align*}
		f_{l}(y)-f_{l}(z)  &  =\liminf_{n\rightarrow\infty}\,f_{n}(y_{n}
		)-\liminf_{n\rightarrow\infty}\,f_{n}(z_{n})\,\geq\,\lim_{n\rightarrow\infty
		}\,f_{k_{n}}(y_{k_{n}})-\liminf_{n\rightarrow\infty}\,f_{k_{n}}(z_{k_{n}})\\
		&  \geq\liminf_{n\rightarrow\infty}\big(f_{k_{n}}(y_{k_{n}})-f_{k_{n}
		}(z_{k_{n}})\big)\geq\liminf_{n\rightarrow\infty}\left\{  -s_{f_{k_{n}}
		}(z_{k_{n}})\,\Vert y_{k_{n}}-z_{k_{n}}\Vert\right\}  \,=\,-s_{f}(z)\,\Vert
		y-z\Vert.
	\end{align*}
	Since $s_{f}$ is locally bounded on $\ri(\dom s_{f})=\ri(\dom\,f_{l})$, we get
	that $f_{l}$ is locally Lipschitz on $\ri(\dom\,f_{l})$. This finishes the
	proof. \hfill$\Box$
	
	\bigskip
	
	\begin{remark}
		Observe that, under the same assumptions, the proof of
		Proposition~\ref{prop: loc lipschitz f_l} shows that that $s_{f_{l}}\leq
		s_{f}$ on $\mathrm{dom}\,s_{f}$.
	\end{remark}
	
	\bigskip
	
	We finish this subsection with the next proposition showing that, by taking a suitable subsequence, we can eliminate the gap between lower and upper epigraphical limits.
	
	\begin{proposition}
		\label{prop: nice subsequence} Let $f,\{f_{n}\}_{n}:\mathbb{R}^{d}
		\rightarrow\mathbb{R}\cup\{+\infty\}$ be proper convex lower semicontinuous
		functions. Assume that $\{s_{f_{n}}\}_{n}$ epigraphically converges to $s_{f}$
		and that there exists some sequence $(x_{n},x_{n}^{\ast},f_{n}(x_{n}
		))\in\triangle f_{n}$ that converges to $(\bar{x},\bar{x}^{\ast},\alpha
		)\in\mathbb{R}^{d}\times\mathbb{R}^{d}\times\mathbb{R}$. Then, 
		$$ f_{l}(\bar x)=f_{u}(\bar x)=\alpha,$$
		and for some strictly increasing sequence $\{k_{n}\}_{n}$ we have
		\[
		f_{l,k_{n}}=f_{u,k_{n}}.
		\]
		
	\end{proposition}
	
	\noindent \textit{Proof.} Let us first observe that thanks to Proposition~\ref{prop: Thibault Combari} we have $f_{l}(x)=f_{u}(x)=\alpha $. Let now 
	\begin{equation*}
		\mathcal{D=\{}z_{i}\}_{i=1}^{\infty }\subset \ri(\dom s_{f})
	\end{equation*}
	be a dense countable set. For $i=1,$ let $\{z_{1,n}\}_{n}\longrightarrow
	z_{1}$ be such that 
	\begin{equation*}
		s_{f_{n}}(z_{1,n})\longrightarrow s_{f}(z_{1})\qquad \text{and}\qquad
		f_{l}(z_{1})=\,\underset{n\rightarrow \infty }{\lim \inf }\,f_{n}(z_{1,n}).
	\end{equation*}
	Take a subsequence $\{k_{1}(n)\}_{n}$ such that:  
	\begin{equation*}
		f_{l}(z_{1})=\,\underset{n\rightarrow \infty }{\lim }
		\,f_{k_{1}(n)}(z_{1,k_{1}(n)})\qquad \text{and}\qquad \partial ^{\circ
		}f_{k_{1}(n)}(z_{1,k_{1}(n)})\longrightarrow \partial ^{\circ
		}f(z_{1}):=z_{1}^{\ast }.
	\end{equation*}
	For $i=2,$ consider a sequence $\{z_{2,n}\}_{n}\longrightarrow z_{2}$ such
	that
	\begin{equation*}
		s_{f_{n}}(z_{2,n})\longrightarrow s_{f}(z_{2}).
	\end{equation*}
	Observe that since $s_{f_{k_1(n)}}\xrightarrow{e} s_f$, Lemma~\ref{lemma: domInclusions} applies and we deduce that $\dom s_f \subset \dom f_{l},_{k_{1}(n)}$. In particular, 
         $f_{l,{k_{1}(n)}}(z_{2})\in\R$.\smallskip\newline
	Replacing $\{z_{2,n}\}_{n}$ by its subsequence $\{z_{2,k_{1}(n)}\}_{n}$ we
	still have $s_{f_{k_{1}(n)}}(z_{i,k_{1}(n)})\longrightarrow s_{f}(z_{i}),$ $i\in \{1,2\}$. Then taking a sub-subsequence $\{k_{2}(k_{1}(n))\}_{n}$ we
	can ensure that  
	\begin{equation*}
		\partial ^{\circ }f_{\left( k_{2}\circ k_{1}\right) (n)}(z_{2,(k_{2}\circ
			k_{1})(n)})\longrightarrow \partial ^{\circ }f(z_{2}):=z_{2}^{\ast }\quad 
		\text{and }\underset{n\rightarrow \infty }{\lim }
		\,f_{(k_{2}\circ k_{1})(n)}(z_{2,(k_{2}\circ k_{1})(n)})\text{ exists in }
		\mathbb{R}.
	\end{equation*}
	We set $\bar{k}_{2}:=k_{2}\circ k_{1}$. Using induction, for every $m>1$, we
	obtain a subsequence $\bar{k}_{m}=k_{m}\circ \ldots \circ k_{1}$ such that
	for all $i\in \{1,\ldots m\}$ we have: 
	\begin{equation*}
		\{z_{i,\bar{k}_{i}(n)}\}_{n}\longrightarrow z_{i}\qquad f_{l}(z_{i})=\,
		\underset{n\rightarrow \infty }{\lim }\,f_{\bar{k}_{i-1}}(z_{i,\bar{k}
			_{i}(n)})\qquad \text{and}\qquad \partial ^{\circ }f_{\bar{k}_{i}(n)}(z_{1,\bar{k}_{i}(n)})\longrightarrow \partial ^{\circ }f(z_{i}):=z_{i}^{\ast }.
	\end{equation*}
	A standard diagonal argument ensures that for every $i\in \mathbb{N}$ the
	sequence $\{\bar{k}_{n}(n)\}_{n\geq i}$ is subsequence of $\{\bar{k}
	_{i}(n)\}_{n\geq i}$. Therefore, thanks to Lemma~\ref{lemma: everySeqTight}
	and the construction, we obtain: 
	\begin{equation*}
		f_{l,\bar{k}_{n}(n)}(z_{i})=\lim_{n\rightarrow \infty }\,f_{\bar{k}
			_{n}(n)}(z_{i,\bar{k}_{n}(n)})=f_{u,\bar{k}_{n}(n)}(z_{i}),\quad \forall
		z_{i}\in \mathcal{D}.
	\end{equation*}
	Since $f_{u,n(k)}$ is convex and lower semicontinuous, using Proposition~\ref{prop: loc lipschitz f_l} and Lemma~\ref{lemma: domInclusions} we deduce
	that 
	\begin{equation*}
		f_{l,\bar{k}_{n}(n)}=f_{u,\bar{k}_{n}(n)}\quad \text{on}\;\ri(\dom s_{f}).
	\end{equation*}
	Thanks to Proposition~\ref{prop: Thibault Combari}, 
	\begin{equation*}
		z_{i}^{\ast }\in \partial f_{u,\bar{k}_{n}(n)}(z_{i})\cap \partial f_{l,\bar{k}_{n}(n)}(z_{i}).
	\end{equation*}
	Let us now define a function $L:\mathbb{R}^{d}\rightarrow \mathbb{R}\cup
	\{+\infty \}$ by 
	\begin{equation*}
		L(z):=\,\sup_{i\in \mathbb{N}}\,\left\{ f_{u,\bar{k}_{n}(n)}(z_{i})+\langle
		z_{i}^{\ast },z-z_{i}\rangle \right\} ,\quad \text{for all }z\in \mathbb{R}^{d}.
	\end{equation*}
	It is straightforward from the definition that $L$ is a lower semicontinuous
	convex function and $L\leq \min \left\{ f_{u,\bar{k}_{n}(n)},~f_{l,\bar{k}_{n}(n)}\right\} $ on the whole space. Notice further that $L=f_{u,\bar{k}_{n}(n)}$ on $\ri(\dom s_{f})$. Then, since $f_{u}$ is convex (by
	Proposition~\ref{Prop: f_u convex}) and $\ri(\dom s_{f})$ is a convex set
	(by Proposition~\ref{prop: convex rint}), we can apply Proposition~\ref{prop: extension to the adherence} to get $L=f_{u,\bar{k}_{n}(n)}$ on $\overline{\dom}(s_{f})$, yielding 
	\begin{equation*}
		f_{l,\bar{k}_{n}(n)}\geq f_{u,\bar{k}_{n}(n)}\quad \text{on}\;\overline{\dom}(s_{f}).
	\end{equation*}
	Finally, thanks to Lemma~\ref{lemma: domInclusions}, we conclude that $f_{l,\bar{k}_{n}(n)}=f_{u,\bar{k}_{n}(n)}$
	on $\mathbb{R}^{d}$. \hfill $\Box $
	

	\section{Main result, final comments and perspectives}
	
	\label{subsec: Ortega}
	
	We are now ready to establish the implications (ii)$\Rightarrow$(i) and
	(iii)$\Rightarrow$(i) of our main result (Theorem~\ref{theo: main}).
	
	\begin{theorem}
		\label{theo: hard implication} Let $f,\{f_{n}\}_{n}:\mathbb{R}^{d}
		\rightarrow\mathbb{R}\cup\{+\infty\}$ be proper convex lower semicontinuous
		functions such that $\inf\,f\in\mathbb{R}$. Assume that $\{s_{f_{n}}\}_{n}$
		epigraphically converges to $s_{f}$ and for some sequence $(x_{n},x_{n}^{\ast
		},f_{n}(x_{n}))\in\triangle f_{n},$ $n\geq1$ we have:
		\[
		\lim_{n\rightarrow\infty}\,(x_{n},x_{n}^{\ast},f_{n}(x_{n}))=(\bar{x},\bar
		{x}^{\ast},f(\bar{x}))\in\triangle f.
		\]
		Then $f_{n}\xrightarrow{e}f$.
	\end{theorem}
	
	\noindent\textit{Proof.} We only need to show that $f\leq f_{l}$ since
	Proposition~\ref{prop: f_u leq f v1} ensures that $f_{u}\leq f$. 
	Let $y\in\dom\,f_{l}$. We
	claim that there exists a sequence $\{y_{n}\}_{n}\subset\mathbb{R}^{d}$ such
	that
	\begin{equation}
		s_{f_{n}}(y_{n})\longrightarrow s_{f}(y)\qquad\text{and}\qquad f_{l}
		(y)=\liminf_{n\rightarrow\infty}\,f_{n}(y_{n}).\label{eq:ortega-3}
	\end{equation}
	Indeed, we distinguish two cases:
	
	\begin{itemize}
		\item \textit{Case 1}\textbf{: }$\;s_{f}(y)<+\infty$.
	\end{itemize}
	
	In this case, since $s_{f_{n}}\xrightarrow{e}s_{f}$, we can choose
	$\{y_{n}\}_{n}$ such that $s_{f}(y)=\lim_{n}s_{f_{n}}(y_{n})$ and apply
	Lemma~\ref{lemma: everySeqTight} to deduce that $f_{l}(y)=\liminf_{n}
	f_{n}(y_{n})$.
	
	\begin{itemize}
		\item \textit{Case 2}: $\;s_{f}(y)=+\infty$.
	\end{itemize}
	
	In this case, every sequence $\{y_{n}\}$ that converges to $y$ should verify
	that $\underset{n\rightarrow\infty}{\lim}s_{f_{n}}(y_{n})=+\infty$. Among
	these sequences, we chose one such that $f_{l}(y)=\underset{n\rightarrow
		\infty}{\liminf}f_{n}(y_{n})$ (\textit{c.f.} Remark~\ref{rem: infsAttained}).
	Therefore (\ref{eq:ortega-3}) holds and the claim is proved.\medskip
	
	Let now $\{k_{1}(n)\}_{n}$ be a strictly increasing subsequence such that
	\[
	f_{l}(y)=\liminf_{n\rightarrow\infty}\,f_{n}(y_{n})\,=\,\lim_{n\rightarrow
		\infty}\,f_{k_{1}(n)}(y_{k_{1}(n)}).
	\]
	Applying Proposition~\ref{prop: nice subsequence} to the sequence
	$\{f_{k_{1}(n)}\}_{n\geq1}$ we get a subsequence $k=k_{2}\circ k_{1}$ of
	$\{k_{1}(n)\}_{n}$, such that $f_{u,k_{n}}=f_{l,k_{n}}=:g$. Observe that
	$f_{k_{n}}$ epi-converges to $g$ and that $g$ is proper convex lower
	semicontinuous. Thus, Theorem~\ref{theo: easy implication} ensures that
	$\{s_{f_{k(n)}}\}_{n}$ epigraphically converges to $s_{g}$. Therefore
	$s_{g}=s_{f}$. Applying Proposition~\ref{prop: f_u leq f v1} to $\{f_{k_{n}
	}\}_{n}$ and $f$, we deduce that $\inf g=\inf f_{u,k_{n}}=\inf\,f\in
	\mathbb{R}$. Thus, we can apply \cite[Corollary 3.1]{PSV2021} (or apply twice
	Theorem~\ref{theo: comparisonPrinciple}) to deduce that $f=g$. In particular,
	\[
	f(y)=g(y)=f_{l,k_{n}}(y)\,\leq\,\lim_{n\rightarrow\infty}\,f_{k_{n}}(y_{k_{n}
	})=f_{l}(y).
	\]
	Since $y$ is arbitrary, we deduce $f\leq f_{l}$. The proof is complete.
	\hfill$\Box$
	
	\bigskip
	
	\begin{theorem}
		\label{theo: without argmin} Let $f,\{f_{n}\}_{n}:\mathbb{R}^{d}
		\rightarrow\mathbb{R}\cup\{+\infty\}$ be proper convex lower semicontinuous
		functions with $\inf\,f\in\mathbb{R}$. Assume that $\{s_{f_{n}}\}_{n}$
		epigraphically converges to $s_{f}$ and
		\[
		\inf\,f_{u}=\inf\,f=\inf\,f_{l}\in\mathbb{R}.
		\]
		Then $\{f_{n}\}_{n}$ epigraphically converges to $f$.
	\end{theorem}
	
	\noindent\textit{Proof.} We follow the arguments of the proof of
	Theorem~\ref{theo: hard implication} with slight modifications. Indeed, it
	suffices to show $f\leq f_{l}$ since Proposition~\ref{prop: fu leq f v2}
	ensures that $f_{u}\leq f$. To this end, we only need to show that $f(y)\leq
	f_{l}(y)$ for all $y\in\dom\,f_{l}$. Fix such $y\in\dom\,f_{l}$ and choose
	again a sequence $\{y_{n}\}_{n}\subset\mathbb{R}^{d}$ such that $s_{f_{n}
	}(y_{n})\rightarrow s_{f}(y)$ and $f_{l}(y)=\liminf_{n}f_{n}(y_{n})$. Take
	$\{k_{n}\}_{n}$ be an increasing subsequence such that
	\[
	f_{l}(y)=\liminf_{n\rightarrow\infty}\,f_{n}(y_{n})=\lim_{n\rightarrow\infty
	}\,f_{k_{n}}(y_{k_{n}}).
	\]
	The main difference, with respect to the proof of
	Theorem~\ref{theo: hard implication}, is that in order to apply
	Proposition~\ref{prop: nice subsequence} to $\{f_{k_{n}}\}_{n}$ we need to
	ensure the existence of a sequence $(x_{k_{n}},x_{k_{n}}^{\ast},f_{k_{n}
	}(x_{k_{n}}))\in\triangle f_{k_{n}}$ that converges (up to a subsequence) to
	some point $(x,x^{\ast},\alpha)$. Since $f$ is proper, there exists at least
	one point $x\in\dom s_{f}$ and sequence $\{x_{m}\}_{m}$ converging to $x$ such
	that $s_{f_{k_{m}}}(x_{k_{m}})$ converges to $s_{f}(x)$. Take $x_{m}^{\ast
	}:=\partial^{\circ}f_{k(m)}(x_{m})$. Using Lemma~\ref{lemma: everySeqTight},
	the hypothesis that $\inf f_{l}=\inf f$, and the fact that $f_{u}\leq f$, we
	can write
	\begin{align*}
		\inf\,f\,\leq\,f_{l}(x)\,\leq\,f_{l,k(m)}(x)  &  =\,\liminf_{m\rightarrow
			\infty}\,f_{k(m)}(x_{m})\\
		&  \leq\,\limsup_{m\rightarrow\infty}\,f_{k(m)}(x_{m})=f_{u,k(m)}
		(x)\,\leq\,\,f_{u}(x)\,\leq\,f(x).
	\end{align*}
	Thus, $\{f_{k(m)}(x_{m})\}_{m}$ is a bounded sequence. We deduce that
	$\{(x_{m},f_{k(m)}(x_{m})),x_{m}^{\ast})\}_{m}$ is also a bounded sequence,
	thus it converges, up to a second subsequence. Therefore, we can apply
	Proposition~\ref{prop: nice subsequence} to the sequence of functions
	$\{f_{k(m)}\}_{m}$. The rest of the proof follows exactly the lines of
	Theorem~\ref{theo: hard implication}. \hfill$\Box$
	

	\begin{remark}
		Due to the fact that in our main result, Theorem~\ref{theo: main}, the limit
		function $f$ is bounded from below, we can slightly generalize it by replacing
		(NC) with the following weaker condition: \smallskip\newline
		\phantom{jo}$\widetilde{(NC)}$\quad There exist $x\in\dom\partial f$ and a sequence $\{x_{n}\}_{n} \subset\mathbb{R}^{d}$ such that: 
		\[
		\lim_{n\to+\infty}(x_{n},f_{n}(x_{n}))=(x,f(x)) \qquad \text{and} \quad \{s_{f_{n}}(x_{n})\}_{n} \quad \text{is bounded.}
		\]
	\end{remark}
	
	\textit{Open problems:} This work is motivated by the celebrated Attouch
	theorem (Theorem~\ref{theo: attouch}), the determination result of slopes 
	\cite{PSV2021}, and the sensitivity result of \cite{DD2023}. All of these
	results are valid in Hilbert spaces, while the first two are also valid in
	Banach spaces (see \cite{AB1993,CT1998} and \cite{TZ2023}). Therefore, a
	natural question is whether Theorem~\ref{theo: main} (our main result) is
	true in Hilbert spaces, or more generally, in reflexive Banach spaces (or
	even in general Banach spaces). While there is no obvious obstruction for
	this extension, the present work relies heavily on local compactness of the
	space, for many of its intermediate results and consequently any potential
	extension should rather rely in a completely different approach.\smallskip 
	\newline
	A much more ambitious project would be to extend the result to pure metric
	spaces, without vector structure, according to the spirit of the
	determination results \cite{DS2022,DLS2023}. One might focus on the notions
	of convexity that have been coined for metric spaces (see, e.g., \cite{AGS-2008}). This is a more challenging task, but the perspective of
	obtaining a metric version of Attouch theorem with its insight on
	variational deviations is tempting and should be explored in the future.


\begin{thebibliography}{99}
		
		\bibitem{AGS-2008}
		L.~Ambrosio, N.~Gigli, and G.~Savar\'{e}.
		\newblock {\em Gradient flows in metric spaces and in the space of probability
			measures}.
		\newblock Lectures in Mathematics ETH Z\"{u}rich. Birkh\"{a}user Verlag, Basel,
		second edition, 2008.
		
		\bibitem{A1977}
		H.~Attouch.
		\newblock Convergence de fonctions convexes, des sous-diff\'{e}rentiels et
		semi-groupes associ\'{e}s.
		\newblock {\em C. R. Acad. Sci. Paris S\'{e}r. A-B}, 284(10):A539--A542, 1977.
		
		\bibitem{A1984-book}
		H.~Attouch.
		\newblock {\em Variational convergence for functions and operators}.
		\newblock Applicable Mathematics Series. Pitman (Advanced Publishing Program),
		Boston, MA, 1984.
		
		\bibitem{AB1993}
		H.~Attouch and G.~Beer.
		\newblock On the convergence of subdifferentials of convex functions.
		\newblock {\em Arch. Math. (Basel)}, 60(4):389--400, 1993.
		
		\bibitem{ABM2014-book}
		H.~Attouch, G.~Buttazzo, and G.~Michaille.
		\newblock {\em Variational analysis in {S}obolev and {BV} spaces}, volume~17 of
		{\em MOS-SIAM Series on Optimization}.
		\newblock Society for Industrial and Applied Mathematics (SIAM), Philadelphia,
		PA; Mathematical Optimization Society, Philadelphia, PA, second edition,
		2014.
		\newblock Applications to PDEs and optimization.
		
		\bibitem{AD-1978}
		H.~Attouch and A.~Damlamian.
		\newblock Strong solutions for parabolic variational inequalities.
		\newblock {\em Nonlinear Anal.}, 2(3):329--353, 1978.
		
		\bibitem{AC-2017}
		D.~Az\'{e} and J.~Corvellec.
		\newblock Nonlinear error bounds via a change of function.
		\newblock {\em J. Optim. Theory Appl.}, 172(1):9--32, 2017.
		
		\bibitem{AC-2004}
		D.~Az\'{e} and J.l Corvellec.
		\newblock Characterizations of error bounds for lower semicontinuous functions
		on metric spaces.
		\newblock {\em ESAIM Control Optim. Calc. Var.}, 10(3):409--425, 2004.
		
		\bibitem{BCD2018}
		T.~Boulmezaoud, P.~Cieutat, and A.~Daniilidis.
		\newblock Gradient flows, second-order gradient systems and convexity.
		\newblock {\em SIAM J. Optim.}, 28(3):2049--2066, 2018.
		
		\bibitem{CT1998}
		C.~Combari and L.~Thibault.
		\newblock On the graph convergence of subdifferentials of convex functions.
		\newblock {\em Proc. Amer. Math. Soc.}, 126(8):2231--2240, 1998.
		
		\bibitem{DDD2015}
		A.~Daniilidis, R.~Deville, E.~Durand-Cartagena, and L.~Rifford.
		\newblock Self-contracted curves in {R}iemannian manifolds.
		\newblock {\em J. Math. Anal. Appl.}, 457(2):1333--1352, 2018.
		
		\bibitem{DD2023}
		A.~Daniilidis and D.~Drusvyatskiy.
		\newblock The slope robustly determines convex functions.
		\newblock {\em Proc. Amer. Math. Soc.}, 151(11):4751--4756, 2023.
		
		\bibitem{DLS2023}
		A.~Daniilidis, T.~M. Le, and D.~Salas.
		\newblock Metric compatibility and determination in complete metric spaces,
		2023.
		\newblock arXiv:2308.14877.
		
		\bibitem{DQ-2022}
		A.~Daniilidis and M.~Quincampoix.
		\newblock Extending {R}ademacher theorem to set-valued maps.
		\newblock {\em hal-03896086}, 2022.
		
		\bibitem{DS2022}
		A.~Daniilidis and D.~Salas.
		\newblock A determination theorem in terms of the metric slope.
		\newblock {\em Proc. Amer. Math. Soc.}, 150(10):4325--4333, 2022.
		
		\bibitem{DgMT-1980}
		E.~De~Giorgi, A.~Marino, and M.~Tosques.
		\newblock Problems of evolution in metric spaces and maximal decreasing curve.
		\newblock {\em Atti Accad. Naz. Lincei Rend. Cl. Sci. Fis. Mat. Nat. (8)},
		68(3):180--187, 1980.
		
		\bibitem{DIL-2015}
		D.~Drusvyatskiy, A.~D. Ioffe, and A.~S. Lewis.
		\newblock Curves of descent.
		\newblock {\em SIAM J. Control Optim.}, 53(1):114--138, 2015.
		
		\bibitem{I-2000}
		A.~D. Ioffe.
		\newblock Metric regularity and subdifferential calculus.
		\newblock {\em Uspekhi Mat. Nauk}, 55(3(333)):103--162, 2000.
		
		\bibitem{L-1988}
		B.~Lemaire.
		\newblock Coupling optimization methods and variational convergence.
		\newblock In {\em Trends in mathematical optimization ({I}rsee, 1986)},
		volume~84 of {\em Internat. Schriftenreihe Numer. Math.}, pages 163--179.
		Birkh\"{a}user, Basel, 1988.
		
		\bibitem{LPT-1995}
		A.~B. Levy, R.~Poliquin, and L.~Thibault.
		\newblock Partial extensions of {A}ttouch's theorem with applications to
		proto-derivatives of subgradient mappings.
		\newblock {\em Trans. Amer. Math. Soc.}, 347(4):1269--1294, 1995.
		
		\bibitem{MP1991}
		P.~Manselli and C.~Pucci.
		\newblock Maximum length of steepest descent curves for quasi-convex functions.
		\newblock {\em Geom. Dedicata}, 38(2):211--227, 1991.
		
		\bibitem{M-1969}
		U.~Mosco.
		\newblock Convergence of convex sets and of solutions of variational
		inequalities.
		\newblock {\em Advances in Math.}, 3:510--585, 1969.
		
		\bibitem{PSV2021}
		P.~P\'{e}rez-Aros, D.~Salas, and E.~Vilches.
		\newblock Determination of convex functions via subgradients of minimal norm.
		\newblock {\em Math. Program.}, 190(1-2):561--583, 2021.
		
		\bibitem{Phelps1993}
		R.~Phelps.
		\newblock {\em Convex functions, monotone operators and differentiability},
		volume 1364 of {\em Lecture Notes in Mathematics}.
		\newblock Springer-Verlag, Berlin, second edition, 1993.
		
		\bibitem{P-1992}
		R.~Poliquin.
		\newblock An extension of {A}ttouch's theorem and its application to
		second-order epi-differentiation of convexly composite functions.
		\newblock {\em Trans. Amer. Math. Soc.}, 332(2):861--874, 1992.
		
		\bibitem{R-1970}
		R.~T. Rockafellar.
		\newblock On the maximal monotonicity of subdifferential mappings.
		\newblock {\em Pacific J. Math.}, 33:209--216, 1970.
		
		\bibitem{R-1990}
		R.~T. Rockafellar.
		\newblock Generalized second derivatives of convex functions and saddle
		functions.
		\newblock {\em Trans. Amer. Math. Soc.}, 322(1):51--77, 1990.
		
		\bibitem{RW1998}
		R.~T. Rockafellar and R.~Wets.
		\newblock {\em Variational {A}nalysis}, volume 317 of {\em Grundlehren der
			mathematischen Wissenschaften [Fundamental Principles of Mathematical
			Sciences]}.
		\newblock Springer-Verlag, Berlin, 1998.
		
		\bibitem{TZ2023}
		L.~Thibault and D.~Zagrodny.
		\newblock Determining functions by slopes.
		\newblock {\em Commun. Contemp. Math.}, 25(7):Paper No. 2250014, 30, 2023.
		
	\end{thebibliography}

	\newpage
	
	\rule{5 cm}{0.5 mm} \bigskip\newline\noindent Aris DANIILIDIS, Sebasti{\'a}n TAPIA-GARC{\'I}A
	
	\medskip
	
	\noindent Institute of Statistics and Mathematical Methods in Economics,
	E105-04 \newline TU Wien, Wiedner Hauptstra{\ss}e 8, A-1040
	Wien\smallskip\newline\noindent E-mail: \{\texttt{aris.daniilidis,
		sebastian.tapia\}@tuwien.ac.at}\newline\noindent
	\texttt{https://www.arisdaniilidis.at/}\newline
	\texttt{https://sites.google.com/view/sebastian-tapia-garcia}
	
	\medskip
	
	\noindent Research supported by the Austrian Science Fund grant \textsc{FWF
		P-36344N}.\newline\vspace{0.2cm}
	
	\noindent David SALAS
	
	\medskip
	
	\noindent Instituto de Ciencias de la Ingenier\'{i}a, Universidad de
	O'Higgins\newline Av. Libertador Bernardo O'Higgins 611, Rancagua, Chile
	\smallskip
	
	\noindent E-mail: \texttt{david.salas@uoh.cl} \newline\noindent
	\texttt{http://davidsalasvidela.cl} \medskip
	
	\noindent Research supported by the grant: \smallskip\newline CMM 
	FB210005 BASAL, FONDECYT 3190229 (Chile)

\end{document}